\documentclass[12pt,oneside,reqno]{amsart}

\usepackage[T1]{fontenc}
\usepackage{amsmath,amsthm}
\numberwithin{equation}{section}
\usepackage{newtxtext,newtxmath}
\usepackage[mathcal]{euscript}
\usepackage{microtype}
\usepackage[a4paper,textwidth=17cm,top=2.5cm,bottom=2.8cm]{geometry}
\usepackage{xcolor}
\definecolor{linkblue}{HTML}{294F70}
\usepackage[
  colorlinks=true,
  linkcolor=blue,
  citecolor=blue,
  urlcolor=blue
]{hyperref}

\theoremstyle{plain}
\newtheorem{theorem}{Theorem}[section]
\newtheorem{proposition}[theorem]{Proposition}
\newtheorem{lemma}[theorem]{Lemma}
\newtheorem{corollary}[theorem]{Corollary}
\newtheorem{conjecture}[theorem]{Conjecture}

\theoremstyle{definition}
\newtheorem{definition}[theorem]{Definition}
\newtheorem{example}[theorem]{Example}

\theoremstyle{remark}
\newtheorem{remark}[theorem]{Remark}

\DeclareMathOperator{\CH}{CH}
\DeclareMathOperator{\cl}{cl}
\DeclareMathOperator{\Res}{Res}
\DeclareMathOperator{\Sing}{Sing}
\DeclareMathOperator{\codim}{codim}

\DeclareMathOperator{\Pic}{Pic}
\DeclareMathOperator{\BB}{BB}

\newcommand{\PP}{\mathbb P}
\newcommand{\alg}{\mathrm{alg}}
\newcommand{\C}{\mathbb C}
\newcommand{\Q}{\mathbb Q}
\newcommand{\F}{\mathcal F}

\newcommand{\Supp}{\operatorname{Supp}}

\title{Baum--Bott Residues and Rational Equivalence}
\author[P. P. Barbassa]{Pedro P. Barbassa}
\address[P. P. Barbassa]{IMECC - UNICAMP\\
Departamento de Matem\'atica\\
Rua S\'ergio Buarque de Holanda, 651\\
13083-970 Campinas-SP, Brazil}
\email{p185839@dac.unicamp.br}

\date{}

\subjclass[2020]{Primary 32S65; Secondary 14C15, 14C17, 14C25}

\keywords{Baum--Bott residues, holomorphic foliations, Chow groups,
rational equivalence, closed rational one-forms}

\begin{document}

\begin{abstract}
We study the Baum--Bott residue formula in the Chow ring of a smooth
complex projective variety. Using the Nash modification and localized
Chern classes, we prove that every Baum--Bott residue of a
codimension-one foliation is represented by an algebraic cycle with
complex coefficients on the singular locus. For foliations defined by
closed rational one-forms, we compute the discrepancy between the
square of the first Chern class of the normal sheaf and the Baum--Bott
cycle. This discrepancy is algebraically trivial and vanishes modulo
rational equivalence for every such form precisely when the variety
has irregularity zero.
\end{abstract}

\maketitle
\tableofcontents

\section{Introduction}

Let $X$ be a smooth connected complex variety of dimension $n\geq 2$.
A foliation $\mathcal F$ of codimension $q$ on $X$ is a coherent
subsheaf $T_{\mathcal F}\subset T_X$ of rank $n-q$ that is closed under
the Lie bracket and saturated in $T_X$. Thus its normal sheaf
$N_{\mathcal F}=T_X/T_{\mathcal F}$ is torsion-free. The locus where
$N_{\mathcal F}$ fails to be locally free is the singular locus
$S=\operatorname{Sing}(\mathcal F)$, which has codimension at least two.

The Baum--Bott theorem expresses characteristic classes of the normal
sheaf in terms of the behavior of the foliation near $S$.
Suppose from now on that $X$ is projective. Let $\varphi$ be a
homogeneous polynomial in the Chern classes, with complex coefficients
and degree $q<d\leq n$, where $c_j$ has degree $j$. To each
connected component $S_\alpha$ of $S$, Baum and Bott associate a class
$\operatorname{Res}_{\varphi}(\mathcal F,S_\alpha)
\in H_{2n-2d}(S_\alpha;\mathbb C)$ that depends only on the germ of
$\mathcal F$ along $S_\alpha$. Their residue formula~\cite{BaumBott1972} reads
\[
  \varphi(N_{\mathcal F})
  =\sum_\alpha\mu_*
    \operatorname{Res}_{\varphi}(\mathcal F,S_\alpha)
  \ \ \text{in }H^{2d}(X;\mathbb C),
\]
where $\mu_{*} \colon H_{2n-2d}(S_{\alpha};\C) \to H_{2n-2d}(X;\C) \to H^{2d}(X;\C)$ and the second arrow is the Poincaré duality.
Subsequent work has clarified how these local classes can be computed
and how they behave under changes of the foliation. Suwa~\cite{Suwa1984} proved the rationality of certain Baum--Bott residues for foliations of complete intersection type. Bracci and Suwa~\cite{BracciSuwa2015} studied the variation of the residues in smooth families.
Corr\^ea and Louren\c{c}o~\cite{CorreaLourenco2019} removed the
genericity assumption in the transverse residue formula of Baum and
Bott~\cite{BaumBott1972}, which expresses the residue in degree $q+1$
as a Grothendieck residue on a transverse $(q+1)$-dimensional ball. Corr\^ea and Suwa~\cite{CorreaSuwaFunctoriality}
established a functoriality formula for generically transverse
holomorphic maps satisfying a compatibility condition on the tangent
and normal sheaves. For a survey of residue theory and its applications to
foliations, see Corr\^ea and Seade \cite{CorreaSeadeIndices2026} and Louren\c{c}o and Reis~\cite{LourencoReis2023}.

On a projective variety one can also ask whether the Baum--Bott formula
holds modulo rational equivalence. The first question concerns the
residue on $S_\alpha$ itself. The cycle class map
\[
  \operatorname{cl}_{S_\alpha}\colon
  \mathrm{CH}_{n-d}(S_\alpha)_{\mathbb C}
  \longrightarrow H_{2n-2d}(S_\alpha;\mathbb C)
\]
sends an algebraic cycle to its homology class. Here the source is the
Chow group of $(n-d)$-dimensional cycles with complex coefficients.
We call the residue \emph{algebraic} if it belongs to the image of this
map, and call any preimage a \emph{Chow representation} of the residue.
Our first result establishes algebraicity in codimension one.

\medskip
\noindent\textbf{Theorem A.}
Let $\mathcal F$ be a codimension-one foliation on $X$. For every
homogeneous characteristic polynomial $\varphi$ with complex coefficients
of degree $2\leq d\leq n$, the Baum--Bott residue
$\operatorname{Res}_{\varphi}(\mathcal F,S_\alpha)$ is algebraic on
each connected component $S_\alpha$ of the singular locus.
\medskip

We obtain Theorem~A as Corollary~\ref{cor:codim-one-all-polynomials}. The proof builds on the use of the
Nash modification by Sert\"oz~\cite{Sertoz1989} and on the comparison of
Nash and Baum--Bott residues by Brasselet and Suwa~\cite{BrasseletSuwa2000}. Localized Chern classes give a
Chow representative for the correction term in this comparison. In
arbitrary codimension $q$, this yields representatives satisfying the
global Chow identity whenever $\varphi$ belongs to the ideal
$(c_{q+1},c_{q+2},\ldots)$. In codimension one the remaining monomials
are the powers of $c_1$. Their residues are obtained from the residue
of $c_1^2$ by taking cap products with powers of the first Chern class
of $\det N_{\mathcal F}$ restricted to the singular locus.

Once the residues have Chow representations, one can ask whether the sum
of their pushforwards recovers the characteristic class in the Chow ring.
For a codimension-one foliation, the case $\varphi=c_1^2$ gives a
comparison independent of any choice of representatives.
Indeed, in dimension $n-2$ the Chow group of the singular locus is
freely generated over $\mathbb C$ by its irreducible components $Z$ of
codimension two in $X$. The residue determines the coefficient of each
such component, denoted by $\operatorname{BB}_{c_1^2}(\mathcal F,Z)$.
We define the \emph{Chow discrepancy} by
\[
  \Delta(\mathcal F)
  :=c_1(N_{\mathcal F})^2
    -\sum_Z\operatorname{BB}_{c_1^2}(\mathcal F,Z)\,[Z]
  \ \ \text{in }\mathrm{CH}^2(X)_{\mathbb C}.
\]
The Baum--Bott theorem gives $\operatorname{cl}_X(\Delta(\mathcal F))=0$.
We compute this discrepancy for foliations defined by closed rational
one-forms. For background on closed rational one-forms, see Pereira~\cite{PereiraClosedForms}.

Let $\eta$ be a nonzero closed rational one-form on $X$, and let
$\mathcal F_\eta$ be the foliation whose tangent sheaf is the saturation
of its kernel in $T_X$. Write $(\eta)_\infty=\sum_i p_iD_i$
for its polar divisor and $\lambda_i=\operatorname{Res}_{D_i}(\eta)$
for its residue along $D_i$. Consider the divisor classes
\[
  \operatorname{Res}(\eta)=\sum_i\lambda_i[D_i],
  \ \ \textnormal{and} \ \
  B_\eta=\sum_{p_i=1}\lambda_i^{-1}[D_i]
  \ \ \text{in }\mathrm{CH}^1(X)_{\mathbb C}.
\]
The second sum runs over the simple poles of $\eta$, along which
closedness implies that the residue is nonzero.

\medskip
\noindent\textbf{Theorem B.}
For every nonzero closed rational one-form $\eta$ on $X$,
\[
  \Delta(\mathcal F_\eta)=\operatorname{Res}(\eta)\cdot B_\eta
  \ \ \text{in }\mathrm{CH}^2(X)_{\mathbb C}.
\]
In particular, $\Delta(\mathcal F_\eta)$ is algebraically
equivalent to zero.
\medskip

The algebraic triviality follows from the fact that the residue divisor
of a closed rational one-form has class in $\operatorname{Pic}^0(X)_{\mathbb C}$.
The discrepancy can nevertheless be nonzero in the Chow group, as we
show on the product of an elliptic curve with the projective line.
More generally, as $\eta$ varies, the classes $\Delta(\mathcal F_\eta)$
span precisely
$\operatorname{Pic}^0(X)_{\mathbb C}\cdot\mathrm{CH}^1(X)_{\mathbb C}$,
the subspace generated by products of an algebraically trivial divisor
class with an arbitrary divisor class. This gives a criterion in terms
of the irregularity $q(X)=h^1(X,\mathcal O_X)$.

\medskip
\noindent\textbf{Corollary C.}
The following conditions are equivalent:
\begin{enumerate}
\item $q(X)=h^1(X,\mathcal{O}_X)=0$;

\item $\Delta(\mathcal{F}_\eta)=0$ for every nonzero
closed rational one-form $\eta$ on $X$.
\end{enumerate}
\medskip

In Section~2 we recall the Baum--Bott residue formula and explain how
the coefficient along a codimension-two component of the singular locus
can be computed on a transverse surface. We also recall the compatibility
of localized Chern classes with the cycle class map. In Section~3 we use
this compatibility to construct a Chow representative for the correction
term relating Nash residues to Baum--Bott residues. In codimension one,
we complete the proof of Theorem~A by expressing the residues for higher
powers of the first Chern class in terms of the residue for its square.
Section~4 studies the Chow discrepancy for foliations defined by closed
rational one-forms. We pass to a birational model on which the form
admits local normal forms and compute separately the contributions from
its zeros and poles. Pushing the resulting identity forward to $X$ gives
Theorem~B. We determine the subspace spanned by all discrepancies
on $X$ arising from closed rational one-forms. This yields the
irregularity criterion of Corollary~C.

\vspace{0.2cm}

\textbf{Acknowledgements.} The author is grateful to Alan Muniz for the many discussions
and for encouraging him to pursue the broader
question of the algebraicity of Baum--Bott residues. He also thanks
Maur\'icio Corr\^ea, Marcos Jardim, and Gabriel Fazoli for helpful
discussions on this work. He is supported by the Coordenação de Aperfeiçoamento de Pessoal de Nível Superior – Brasil (CAPES) – Finance Code 001 and the PDSE Program No. 88881.220396/2025-01.
\vspace{0.2cm}

\textbf{AI disclosure.} During the development and preparation of this article, the author used ChatGPT-6 and 6.2 Astra as an interactive tool to explore proof strategies, improve the writing, and search for references. The author assumes full responsibility for the content of this article.

\section{Foliations, Baum--Bott residues and localization of classes}

The main purpose of this section is to fix the notation and conventions used through this text. We will always use $X$ to denote a smooth connected projective variety over $\C$ of dimension $n \geq 2$.

\subsection{Foliations and Baum--Bott residues}
\begin{definition}\label{def: foliation}
A codimension $q$ foliation $\mathcal F$ on $X$
is a coherent subsheaf $T_{\mathcal F}\subset T_X$ of rank $n-q$
such that $T_X/T_{\mathcal F}$ is torsion-free and $T_{\F}$ is closed under the Lie bracket
\(
    [T_{\mathcal F},T_{\mathcal F}]\subset T_{\mathcal F}.
\)
The first condition says that $T_{\mathcal F}$ is saturated and
the second is the integrability condition.
\end{definition}

We define the singularities of the foliation $S =\Sing(\F)$ as the points of $X$ such that the quotient $T_X/T_\F$ fails to be locally free. Since the quotient is torsion-free, we have $\codim (\Sing(\F)) \geq 2$. The foliation is called regular if has no singularities.
Though this text we use the notation $N_\F : = T_X/T_\F$.

Notice that the foliation corresponds to a section $\omega\in H^0\big(X,\Omega_X^q\otimes \det(N_{\mathcal F})\big)$ satisfying that at a generic point of $X$, there exist germs of one-forms $\omega_i$ with $1 \leq i \leq q$, such that
\[
\omega = \omega_1 \wedge \dots \wedge \omega_q \ \ \textnormal{and} \ \ d\omega_i \wedge \omega = 0 \ \ \textnormal{for every} \  i =1,\dots,q. 
\]

\begin{example}
Let $X=\mathbb P^1\times\mathbb P^1$ and consider the diagonal
$\mathbb C^*$-action given in affine coordinates by
$t\cdot(x,y)=(tx,ty)$. Its infinitesimal generator is the global
vector field
\[
  v=x\frac{\partial}{\partial x}+y\frac{\partial}{\partial y}.
\]
The orbits define a foliation $\mathcal F$ with rational first
integral $y/x$. On the dense torus $(\mathbb C^*)^2$, it is defined
by the closed logarithmic one-form
\[
  \eta=\frac{dy}{y}-\frac{dx}{x}.
\]
The zero scheme of $v$ consists of the four reduced points
$Z=\{0,\infty\}\times\{0,\infty\}$. Taking the wedge product with
$v$ identifies the quotient $T_X/\mathcal O_X v$ with
$\mathcal I_Z\otimes\det T_X$. Thus
\[
  T_{\mathcal F}\simeq\mathcal O_X \ \ \textnormal{and} \ \
  N_{\mathcal F}\simeq\mathcal I_Z\otimes\mathcal O_X(2,2).
\]
\end{example}

\begin{definition}
Let $f\colon Y\to X$ be a dominant morphism between smooth complex
varieties, and let $\mathcal F$ be a codimension $q$ foliation on $X$.
The \emph{pullback foliation} $f^*\mathcal F$ is the codimension $q$ foliation on $Y$
whose tangent sheaf is the saturation in $T_Y$ of
\[
    \ker\!\left(
        T_Y \xrightarrow{\,df\,} f^*T_X
            \longrightarrow f^*N_{\mathcal F}
    \right).
\]
\end{definition}
If $\mathcal F$ is defined by a twisted $q$-form
$\omega\in H^0(X,\Omega_X^q\otimes L)$, then its pullback
$f^*\omega\in H^0(Y,\Omega_Y^q\otimes f^*L)$ defines
$f^*\mathcal F$ on a dense open subset of $Y$. Writing $D$ for
the divisorial zero part of $f^*\omega$, we obtain a global
defining form for $f^*\mathcal F$ by viewing $f^*\omega$ as a
section of $\Omega_Y^q\otimes f^*L(-D)$.
\vspace{0.2cm}

On the regular points of the foliation $U := X \setminus S$, the normal sheaf carries the Bott partial connection along the leaves:
\[
    \nabla^{\mathrm B}_v(\overline w)
    =\overline{[v,w]},
\]
for $v\in T_{\mathcal F}|_U$ and $w$ is a local lift to $T_U$ of a section
$\overline w\in N_{\mathcal F}|_U$; for partial connections see \cite{FazoliJetsI,FazoliJetsII}.
Integrability makes this well-defined and flat along the leaves.
Let $\varphi$ be a homogeneous polynomial of
degree $d$ with complex coefficients in the Chern
classes $c_j$, with $c_j$ having degree $j$. Assuming $q < d \leq n$, Bott vanishing gives
\[
    \varphi(N_{\mathcal F}|_U)=0
    \quad\text{in }H^{2d}(U,\mathbb C).
\]
Using the vanishing on the regular locus, the Baum--Bott construction
produces a localized class in $H^{2d}(X,U;\mathbb C)$ whose image in
$H^{2d}(X;\mathbb C)$ is $\varphi(N_{\mathcal F})$. Decomposing this
class according to the connected components of the singular locus
gives the following.

\begin{theorem}[Baum--Bott {\cite{BaumBott1972}}]
\label{thm:classical-baum-bott-c1-square}
For each connected component $S_{\alpha}\subset S$, there are canonical classes 
$\Res_{\varphi}(\F,S_{\alpha}) \in H_{2n-2d}(S_{\alpha};\C)$ depending only on $\varphi$ and the germ of $\F$ along $S_{\alpha}$, such that
\begin{equation}\label{eq:classical-baum-bott}
    \varphi(N_{\mathcal F})
    =
    \sum_{\alpha}
    \mu_*\Res_{\varphi}(\mathcal F,S_{\alpha})
    \quad\text{in }H^{2d}(X,\mathbb C),
\end{equation}
where $\mu_{*} \colon H_{2n-2d}(S_{\alpha};\C) \to H_{2n-2d}(X;\C) \to H^{2d}(X;\C)$ and the later map is the Poincaré duality.
\end{theorem}

The elements $\Res_{\varphi}(\F,S_{\alpha})$ are called the \emph{Baum-Bott residues} of the foliation $\F$.
Suppose now that $\mathcal F$ has codimension one and take
$\varphi=c_1^2$. For each connected component $S_\alpha$ of
$\operatorname{Sing}(\mathcal F)$, the fundamental classes of the
irreducible components $Z\subseteq S_\alpha$ of codimension two in $X$
form a basis of $H_{2n-4}(S_\alpha;\mathbb C)$. Thus the Baum--Bott
residue has a unique expansion
\begin{equation}\label{eq: BB coef}
  \operatorname{Res}_{c_1^2}(\mathcal F,S_\alpha)
  =
  \sum_{\substack{Z\ \mathrm{irred.\ component\ of}\ S_\alpha\\
                  \operatorname{codim}_X Z=2}}
  \operatorname{BB}_{c_1^2}(\mathcal F,Z)\,[Z]
  \quad\text{in } H_{2n-4}(S_\alpha;\mathbb C),
\end{equation}
where $[Z]$ denotes the image of the fundamental class of $Z$ in
$S_\alpha$. We call the scalar
$\operatorname{BB}_{c_1^2}(\mathcal F,Z)\in\mathbb C$ the
\emph{Baum--Bott coefficient} of $\mathcal F$ along $Z$.
If $\operatorname{codim}_X S_\alpha>2$, the homology group above
vanishes and the residue is zero.

\subsection{Grothendieck residue}
We recall the local Grothendieck residue and use it to express the
Baum--Bott coefficient of a codimension-two singular component as
a residue on a transverse surface.

Let $f_1,\ldots,f_m,h$ be holomorphic germs at $0\in\mathbb C^m$,
and suppose that $0$ is an isolated common zero of
$f_1,\ldots,f_m$. In local coordinates $z_1,\ldots,z_m$, the
Grothendieck residue with respect to the ordered tuple
$(f_1,\ldots,f_m)$ is
\[
\operatorname{Res}_0\!
\left[
\frac{h\,dz_1\wedge\cdots\wedge dz_m}{f_1,\ldots,f_m}
\right]
:=
\frac{1}{(2\pi i)^m}
\int_{\Gamma_\varepsilon}
\frac{h\,dz_1\wedge\cdots\wedge dz_m}{f_1\cdots f_m}.
\]
Here $\Gamma_\varepsilon=\{|f_j|=\varepsilon_j,\ 1\leq j\leq m\}$
is taken in a sufficiently small neighborhood of the origin, for
sufficiently small positive regular values, and is oriented by
$d\arg f_1\wedge\cdots\wedge d\arg f_m>0$. The residue depends
only on the class of $h$ modulo $(f_1,\ldots,f_m)$. If the Jacobian
matrix $J_f(0)=(\partial f_i/\partial z_j)(0)$ is invertible, then
the residue equals $h(0)/\det J_f(0)$. 

Suppose that $\mathcal F$ has codimension one, and let $Z$ be an
irreducible component of $\operatorname{Sing}(\mathcal F)$ of
codimension two in $X$. Choose a general point $p\in Z$ at which
$Z$ is smooth and which lies on no other irreducible component of
the singular locus. Let $i\colon B\hookrightarrow X^{\mathrm{an}}$
be a general, sufficiently small, two-dimensional complex ball
through $p$, transverse to $Z$ at $p$. Restricting a local defining
one-form of $\mathcal F$ to $B$ gives a foliation $\mathcal F_B$
with an isolated singularity at $p$. By the transverse residue
formula of Corr\^ea and Louren\c{c}o~\cite[Theorem 1.2]{CorreaLourenco2019},
\[
  \operatorname{BB}_{c_1^2}(\mathcal F,Z)
  =
  \operatorname{BB}_{c_1^2}(\mathcal F_B,p).
\] 

Let $\mathcal F$ be a germ of a holomorphic foliation on a smooth
complex surface, with an isolated singularity at $0$. In local
coordinates $(x,y)$ centered at $0$, choose a generator
$v=a\,\partial_x+b\,\partial_y$ of $T_{\mathcal F}$ and write
$J_v=\partial(a,b)/\partial(x,y)$. Since $a$ and $b$ have an isolated
common zero at $0$, the Baum--Bott coefficient is given by the
Grothendieck residue~\cite{BaumBott1972}
\begin{equation}\label{eq:bb-surface-grothendieck}
  \mathrm{BB}_{c_1^2}(\mathcal F,0)
  =
  \operatorname{Res}_0\!
  \left[
    \frac{(\operatorname{tr}J_v)^2\,dx\wedge dy}{a,b}
  \right].
\end{equation}
This formula holds for every isolated singularity, including a
degenerate one. If the singularity is nondegenerate, the linearization
$J_v(0)$ has nonzero eigenvalues $\alpha$ and $\beta$, and the formula
simplifies to
\[
  \mathrm{BB}_{c_1^2}(\mathcal F,0)
  =
  \frac{(\operatorname{tr}J_v(0))^2}{\det J_v(0)}
  =
  \frac{(\alpha+\beta)^2}{\alpha\beta}.
\]
In Subsection~\ref{subsec:residue-calculations-proof}, we use the
eigenvalue formula when the transverse singularity is nondegenerate.
At the intersection of a simple pole with a pole of higher order,
the transverse singularity is degenerate, and we compute its
coefficient directly from~\eqref{eq:bb-surface-grothendieck}.

\subsection{Localized Chern classes and the cycle class map}
For local Chern classes, see~\cite{Iversen1976}.
For a complex algebraic variety $T$, let $\mathrm{CH}_k(T)$ denote the
group of $k$-dimensional cycles modulo rational equivalence, and write
$\mathrm{CH}_k(T)_A=\mathrm{CH}_k(T)\otimes_{\mathbb Z}A$ for
$A=\mathbb Q$ or $\mathbb C$. Since $X$ is smooth of dimension $n$,
cap product with its fundamental class gives the isomorphism
\[
  \mathrm{CH}^a(X)_A \xrightarrow{\ \cap[X]\ } \mathrm{CH}_{n-a}(X)_A.
\]
The cycle class map
\[
  \operatorname{cl}_T\colon \mathrm{CH}_k(T)_A
  \longrightarrow H^{\mathrm{BM}}_{2k}(T^{\mathrm{an}};A)
\]
sends the class of an integral subvariety to the pushforward of its
Borel--Moore fundamental class~\cite[\S~19.1]{Fulton1998}. When $T$
is proper, $T^{\mathrm{an}}$ is compact and Borel--Moore homology
agrees with ordinary singular homology. This applies to $X$ and to
every closed subvariety of $X$. We suppress analytification in the
notation for homology and cohomology. On $X$, Poincar\'e duality gives
the cohomological cycle class map
\[
  \operatorname{cl}_X\colon \mathrm{CH}^a(X)_A\longrightarrow H^{2a}(X;A).
\]

We write $\mathrm{CH}^a(X)_{\mathrm{hom},A}=\ker(\operatorname{cl}_X)$.
The subspace $\mathrm{CH}^a(X)_{\mathrm{alg},A}$ of algebraically
trivial classes is spanned over $A$ by differences
$i_{t_1}^*\Gamma-i_{t_0}^*\Gamma$, where $C$ is a smooth connected
algebraic curve, $\Gamma\in\mathrm{CH}^a(X\times C)$, and
$i_t\colon X\hookrightarrow X\times C$ is the inclusion of the fiber
over $t\in C$. The inclusions $i_{t_0}$ and $i_{t_1}$ are homotopic
after analytification, so
\[
  \mathrm{CH}^a(X)_{\mathrm{alg},A}
  \subseteq \mathrm{CH}^a(X)_{\mathrm{hom},A}.
\]
This inclusion can be strict in codimension two over either coefficient
field. Griffiths exhibited two lines on a smooth quintic threefold in
$\mathbb P^4_{\mathbb C}$ whose difference is homologically trivial and
has infinite order modulo algebraic equivalence~\cite[\S~14]{Griffiths1969II}.

Let $j\colon T\hookrightarrow Y$ be a closed subvariety of a smooth
projective complex variety of dimension $n$, and put $U=Y\setminus T$.
Let $K^\bullet$ be a bounded complex of vector bundles on $Y$, exact
on $U$, with class $e=\sum_i(-1)^i[K^i]\in K^0(Y)$. For each $a\geq 1$,
we construct a Chow class on $T$ whose pushforward is $c_a(e)\cap[Y]$
and identify its homology class on $T$.

Write $H_T^m(Y;A)=H^m(Y,U;A)$. Capping with the fundamental class
of $Y$ in its complex orientation gives the Alexander duality
isomorphism
\[
  \operatorname{AD}_{Y,T}\colon H_T^m(Y;A)
  \xrightarrow{\ \sim\ } H_{2n-m}(T;A).
\]
For $u\in H^{2i}(Y;A)$ and $v\in H_T^{2k}(Y;A)$, the relative cup
product satisfies
\[
  \operatorname{AD}_{Y,T}(u\smile v)
  =j^*u\cap\operatorname{AD}_{Y,T}(v).
\]

Fulton's graph construction associates to $K^\bullet$ its localized
Chern character $\operatorname{ch}_T^Y(K^\bullet)$%
~\cite[\S~18.1]{Fulton1998}. Evaluating its degree-$a$ component on
$[Y]$ gives
\[
  \gamma_a(K^\bullet)
  :=\operatorname{ch}_{T,a}^Y(K^\bullet)\cap[Y]
  \in\mathrm{CH}_{n-a}(T)_{\mathbb Q}.
\]
By \cite[Proposition~18.1(a)]{Fulton1998},
\begin{equation}\label{eq:localized-character-pushforward}
  j_*\gamma_a(K^\bullet)=\operatorname{ch}_a(e)\cap[Y].
\end{equation}
The analytified complex $K^{\bullet,\mathrm{an}}$ determines a relative
topological $K$-class for the pair $(Y^{\mathrm{an}},U^{\mathrm{an}})$
whose image in $K^0(Y^{\mathrm{an}})$ is $e^{\mathrm{an}}$. We write
$\operatorname{ch}^{T,\mathrm{top}}_a(K^{\bullet,\mathrm{an}})$ for
the degree-$2a$ component of its Chern character and, for $a\geq 1$,
$c_a^{T,\mathrm{top}}(K^{\bullet,\mathrm{an}})$ for its $a$th Chern
class, both in $H_T^{2a}(Y;\mathbb Q)$.

\begin{lemma}\label{lem:localized-comparison}
For every $a\geq 0$,
\[
  \operatorname{cl}_T\bigl(\gamma_a(K^\bullet)\bigr)
  =\operatorname{AD}_{Y,T}\!\left(
    \operatorname{ch}^{T,\mathrm{top}}_a(K^{\bullet,\mathrm{an}})
  \right).
\]
\end{lemma}

\begin{proof}
The cycle class map commutes with localized Chern characters
\cite[Example~19.2.6]{Fulton1998}. Applying this compatibility to $[Y]$
in degree $a$ gives the cap product of
$\operatorname{ch}^{T,\mathrm{top}}_a(K^{\bullet,\mathrm{an}})$
with the topological fundamental class of $Y$. This is its Alexander dual.
\end{proof}

We use Newton's identities to localize the Chern classes with rational
coefficients. For the direct graph construction, see
\cite[Example~18.1.3]{Fulton1998}. For $a\geq 1$, set
\begin{equation}\label{eq:localized-chern-definition}
  \Gamma_a(K^\bullet)
  :=\frac{1}{a}\sum_{r=1}^a(-1)^{r-1}r!\,
    c_{a-r}(j^*e)\cap\gamma_r(K^\bullet)
  \in\mathrm{CH}_{n-a}(T)_{\mathbb Q}.
\end{equation}
Here $c_i(j^*e)$ is the Chern class operator of the restricted virtual
bundle on $T$, and $c_0$ acts as the identity.

\begin{corollary}\label{cor:localized-chern-comparison}
For every $a\geq 1$, the class $\Gamma_a(K^\bullet)$ satisfies
\begin{equation}\label{eq:localized-chern-pushforward}
  j_*\Gamma_a(K^\bullet)=c_a(e)\cap[Y]
  \qquad\text{in }\mathrm{CH}_{n-a}(Y)_{\mathbb Q}
\end{equation}
and
\begin{equation}\label{eq:localized-chern-comparison}
  \operatorname{cl}_T\bigl(\Gamma_a(K^\bullet)\bigr)
  =\operatorname{AD}_{Y,T}\!\left(
     c_a^{T,\mathrm{top}}(K^{\bullet,\mathrm{an}})
   \right)
  \qquad\text{in }H_{2n-2a}(T;\mathbb Q).
\end{equation}
\end{corollary}

\begin{proof}
The projection formula and \eqref{eq:localized-character-pushforward}
give
\[
  j_*\Gamma_a(K^\bullet)
  =\left(\frac{1}{a}\sum_{r=1}^a(-1)^{r-1}r!\,
    c_{a-r}(e)\operatorname{ch}_r(e)\right)\cap[Y]
  =c_a(e)\cap[Y],
\]
where the last equality is Newton's identity.

For the comparison on $T$, recall that if $u,v\in H_T^*(Y;\mathbb Q)$
and $\bar u$ is the image of $u$ in $H^*(Y;\mathbb Q)$, then
$u\smile v=\bar u\smile v$ in $H_T^*(Y;\mathbb Q)$. Thus in Newton's
identity for the relative $K$-class of $K^{\bullet,\mathrm{an}}$,
we may forget support in the Chern class factor and retain it in the
Chern character factor. This gives
\[
  a\,c_a^{T,\mathrm{top}}(K^{\bullet,\mathrm{an}})
  =\sum_{r=1}^a(-1)^{r-1}r!\,
    c_{a-r}(e^{\mathrm{an}})\smile
    \operatorname{ch}^{T,\mathrm{top}}_r(K^{\bullet,\mathrm{an}}).
\]
By the cup product formula above and Lemma~\ref{lem:localized-comparison},
Alexander duality takes each cup product on the right to
$j^*c_{a-r}(e^{\mathrm{an}})\cap\operatorname{cl}_T(\gamma_r(K^\bullet))$.
The cycle class map
commutes with Chern class operators~\cite[Proposition~19.1.2]{Fulton1998},
so the sum of these terms with their coefficients is
$a\,\operatorname{cl}_T(\Gamma_a(K^\bullet))$ by
\eqref{eq:localized-chern-definition}. This proves
\eqref{eq:localized-chern-comparison}.
\end{proof}

\section{Algebraicity of Baum–Bott residues}

\subsection{The lifting problem}
Let $\F$ be a codimension $q$ foliation on $X$, and let 
$\varphi$ be a homogeneous characteristic polynomial of degree $q < d \leq n$. For every connected components  $S_{\alpha} \subset S = \Sing(\F)$, consider the cycle class map 
\[
\cl_{S_{\alpha}} \colon \CH_{n-d}(S_{\alpha})_\C \to H_{2n-2d}(S_{\alpha}; \C).
\]

We say that the Baum--Bott residue $\Res_{\varphi}(\F,S_{\alpha})$ is \emph{algebraic} if there exists an element $R_{\varphi,\alpha} \in \CH_{n-d}(S_{\alpha})_{\C}$ such that
\[
\cl_{S_{\alpha}}(R_{\varphi,\alpha}) = \Res_{\varphi}(\F,S_{\alpha}).
\]
We call such a class \(R_{\varphi,\alpha}\) a \emph{Chow representation} of
the residue. Its existence depends on whether the residue lies in
the image of \(\mathrm{cl}_{S_\alpha}\), while its uniqueness, when
it exists, depends on the kernel of this map. Given the residue $\Res_{\varphi}(\F,S_{\alpha}) \in H_{2n-2d}(S_{\alpha}; \C)$, we denote its localized class by
\[
\rho_{\varphi,\alpha} := \operatorname{AD}_{X,S_\alpha}^{-1}(\Res_{\varphi}(\F,S_\alpha)) \in H_{S_\alpha}^{2d}(X;\C).
\]
\begin{lemma}\label{lem: trivial-cases}
    Let $\F$ be a codimension $q$ foliation on $X$, and let $\varphi$ be a homogeneous characteristic polynomial of degree $q < d \leq n$. If
    \begin{enumerate}
        \item $d < \codim_X(S_{\alpha})$, then the Baum--Bott residue is trivial;

        \item $d = \codim_X(S_{\alpha})$, then there exist a unique $R_{\varphi,\alpha}$ such that $\cl_{S_{\alpha}}(R_{\varphi,\alpha}) = \Res_{\varphi}(\F,S_{\alpha})$;

        \item $d =\dim X = n$, then the residue admits a Chow representative.
    \end{enumerate}
\end{lemma}
\begin{proof}
    All the results follow from the cycle class map. If $d < \codim_X(S_{\alpha)}$ then both groups $H_{2n-2d}(S_{\alpha}; \C)$ and $\CH_{n-d}(S_{\alpha})_\C$ are trivial. If $d = \codim_X(S_{\alpha})$, then both the Chow group and the homology group are generated by the classes of the irreducible components of $S_{\alpha}$ with the same dimension as $S_{\alpha}$, therefore the cycle class map is an isomorphism.

    If $d = n$, since $S_{\alpha}$ is a connected projective variety, its has $H_0(S_{\alpha};\C)= \C$. So the cycle class map in this case is the complexified degree map
    \[
    \cl_{S_{\alpha}} \colon \CH_0(S_{\alpha})_{\C} \to \C,
    \]
    which is surjective and therefore follows the result.
\end{proof}

\begin{remark}
    Since in general we cannot guarantee that a Chow representation is unique, given $R_{\varphi,\alpha}$ all other Chow representations of $\Res_{\varphi}(\F,S_{\alpha})$ are of the form $R_{\varphi,\alpha} + K_\alpha$ with $K_\alpha \in \ker(\cl_{S_{\alpha}})$.
\end{remark}

To study the lifting problem beyond the cases covered by Lemma~\ref{lem: trivial-cases},
we use the Nash modification of the foliation. On a resolution of
this modification, we construct a Chow class that represents the
difference between the localized Baum--Bott and Nash classes.

\subsection{The Nash modification and the correction class}
\label{subsec:nash-modification}
Sert\"oz compared Baum--Bott and Nash residues in homology and expressed their difference by means of the Grassmann graph construction \cite[Theorems~IV.4 and~V.1]{Sertoz1989}. Brasselet and Suwa defined Nash residues on the possibly singular Nash modification and established a comparison formula for the
localized characteristic classes in cohomology with support
\cite[\S\S~3--4, especially equation~(4.6)]{BrasseletSuwa2000}.
We use this formula in the projective algebraic case. 

We keep the notation that $\F$ is a codimension $q$ foliation on $X$ and $S = \Sing(\F)$. Over $U = X \setminus S$, the tangent bundle of the regular foliation defines
a section
\[
U\longrightarrow \operatorname{Gr}_{n-q}(T_X)
\ \ \textnormal{given by } \ \
x\longmapsto T_{\mathcal F,x}.
\]
Let $X^\nu$ be the reduced closure of its image, and denote the
resulting projective birational morphism by
\(
\pi^\nu:X^\nu\longrightarrow X.
\)
The restriction of the universal sequence on the Grassmann
bundle gives
\[
0\longrightarrow V^\nu
\longrightarrow (\pi^\nu)^*T_X
\longrightarrow Q^\nu
\longrightarrow 0,
\]
where $V^\nu$ and $Q^\nu$ are vector bundles of ranks $n-q$
and $q$, respectively. Over $(\pi^\nu)^{-1}(U)\simeq U$, the
bundle $Q^\nu$ identifies with $\mathcal N_{\mathcal F}|_U$.
Set
\[
T^\nu=(\pi^\nu)^{-1}(S)_{\mathrm{red}}
\ \ \textnormal{and } \ \
T^\nu_\alpha=(\pi^\nu)^{-1}(S_\alpha)_{\mathrm{red}}.
\]
The variety $X^\nu$ need not be smooth. Nevertheless, the Bott
partial connection on $\mathcal N_{\mathcal F}|_U$ yields
localized characteristic classes of $Q^\nu$. These are the
classes used to define Nash residues
\cite[Section~3]{BrasseletSuwa2000}.
More precisely, for a homogeneous characteristic polynomial
$\varphi$ of degree $q<d\le n$, denote the corresponding
localized class by
\[
\nu^\nu_{\varphi,\alpha}
\in
H_{T^\nu_\alpha}^{2d}(X^\nu;\mathbb C).
\]
Its image under the Alexander homomorphism for $X^\nu$ is the
sum of the Nash residues supported on $T^\nu_\alpha$.
On a singular variety this Alexander homomorphism need not
be an isomorphism.
Choose a projective resolution
\(
r:Y\longrightarrow X^\nu
\)
which is an isomorphism over $U$, and put
\[
\pi=\pi^\nu\circ r:Y\longrightarrow X.
\]
Such a resolution exists in characteristic zero. In particular,
$Y$ is smooth and projective of dimension $n$.
Fixing the notation, we have $T=\pi^{-1}(S)_{\mathrm{red}}$, $T_\alpha=\pi^{-1}(S_\alpha)_{\mathrm{red}}$.
Then
\[
T=\coprod_\alpha T_\alpha
\ \ \textnormal{and } \ \
Y\setminus T\simeq U.
\]
Pulling back the tautological sequence gives
\[
0\longrightarrow V
\longrightarrow \pi^*T_X
\longrightarrow Q
\longrightarrow 0,
\]
where $V=r^*V^\nu$ and $Q=r^*Q^\nu$. Define
\[
\nu_{\varphi,\alpha}
=
r^*\nu^\nu_{\varphi,\alpha}
\in H_{T_\alpha}^{2d}(Y,\mathbb C).
\]
Thus $\nu_{\varphi,\alpha}$ denotes the pulled-back Nash
localized class on the smooth variety $Y$.

We now use the construct of the correction between the derived pullback
of $\mathcal N_{\mathcal F}$ and the tautological quotient $Q$; see \cite[Section 3]{BrasseletSuwa2000}.
Choose a finite locally free resolution
\[
0\longrightarrow E_m
\longrightarrow\cdots
\longrightarrow E_1
\longrightarrow T_{\mathcal F}
\longrightarrow0.
\]
Together with $T_{\mathcal F}\hookrightarrow T_X$, it gives
a locally free resolution of $\mathcal N_{\mathcal F}$.
On $Y$, consider
\[
K^\bullet:
\quad
0\longrightarrow\pi^*E_m
\longrightarrow\cdots
\longrightarrow\pi^*E_1
\longrightarrow\pi^*T_X
\longrightarrow Q
\longrightarrow0,
\]
where $\pi^*T_X$ is placed in degree zero and $Q$ in degree
one. The composite
\(
\pi^*E_1\rightarrow Q
\)
vanishes over $Y\setminus T$ and therefore vanishes everywhere.
Thus $K^\bullet$ is a complex, which is exact on
$Y\setminus T$.
Set
\[
\xi=[L\pi^*\mathcal N_{\mathcal F}]
\ \ \textnormal{and } \ \
e=[K^\bullet].
\]
In $K_0(Y)$ we have
\(
\xi=[Q]+e.
\)
The derived pullback is represented by the pulled-back
locally free resolution.
Let $\varphi$ be any homogeneous characteristic polynomial
of degree $d$, with rational coefficients.
The identity
\(
c(Q+e)=c(Q)c(e)
\)
gives an expansion
\[
\varphi(Q+e)-\varphi(Q)
=
\sum_\lambda
A_\lambda(c(Q))
\prod_{r=1}^{t_\lambda}c_{a_{\lambda,r}}(e),
\]
where $t_\lambda\ge1$ and $a_{\lambda,r}\ge1$.
Here $A_\lambda(c(Q))$ is a polynomial in the Chern
classes of $Q$, homogeneous of degree
\[
d-\sum_{r=1}^{t_\lambda}a_{\lambda,r}.
\]
Every summand contains a positive-degree Chern class of $e$,
because the difference vanishes formally when $e=0$.
Put
\[
P_\lambda
=
A_\lambda(c(Q))
\prod_{r=2}^{t_\lambda}c_{a_{\lambda,r}}(e)
\in \CH^{d-a_{\lambda,1}}(Y)_{\mathbb Q}.
\]
The empty product is understood to be $1$.
The symbol $P_{\lambda,T}$ will denote the characteristic
class operator on $\CH_*(T)_{\mathbb Q}$ obtained by restricting
the vector bundles in this expression to $T$. Thus
\[
P_{\lambda,T}:
\CH_k(T)_{\mathbb Q}
\longrightarrow
\CH_{k-d+a_{\lambda,1}}(T)_{\mathbb Q}.
\]
Define
\[
D_\varphi
:=
\sum_\lambda
P_{\lambda,T}\cap
\Gamma_{a_{\lambda,1}}(K^\bullet)
\in \CH_{n-d}(T)_{\mathbb Q},
\]
where $\Gamma_{a_{\lambda,1}}(K^\bullet)$ is defined in \eqref{eq:localized-chern-definition}.
Write $D_{\varphi,\alpha}$ for the component of $D_\varphi$
in $\CH_{n-d}(T_\alpha)_{\mathbb Q}$.

\begin{proposition}
\label{prop:nash-correction}
For $j\colon T \hookrightarrow Y$ and every homogeneous characteristic polynomial $\varphi$
of degree $d$,
\[
j_*D_\varphi
=
\bigl(\varphi(\xi)-\varphi(Q)\bigr)\cap[Y] \in \CH_{n-d}(Y)_{\Q}.
\]
If moreover $q<d\le n$, then, for every $\alpha$,
\[
\operatorname{cl}_{T_\alpha}(D_{\varphi,\alpha})
=
\operatorname{AD}_{Y,T_\alpha}
\bigl(
\pi^*\rho_{\varphi,\alpha}
-
\nu_{\varphi,\alpha}
\bigr).
\]
\end{proposition}

\begin{proof}
By the projection formula and Corollary~\ref{cor:localized-chern-comparison},
\[
\begin{aligned}
j_*D_\varphi
&= \sum_\lambda P_\lambda \cap
   j_*\Gamma_{a_{\lambda,1}}(K^\bullet) 
= \left(\sum_\lambda P_\lambda\,c_{a_{\lambda,1}}(e)\right)
   \cap [Y] \\
&= \bigl(\varphi(\xi)-\varphi(Q)\bigr)\cap [Y],
\end{aligned}
\]
where the last equality follows from the Whitney expansion and
the identity $\xi=[Q]+e$.
Assume now that $q<d\leq n$. Let $\rho_\varphi\in H^{2d}_S(X;\mathbb C)$
be the class whose component on $S_\alpha$ is
\[
\rho_{\varphi,\alpha}
=
\operatorname{AD}_{X,S_\alpha}^{-1}
\bigl(\operatorname{Res}_\varphi(\mathcal F,S_\alpha)\bigr),
\]
and let $\nu_\varphi\in H^{2d}_T(Y;\mathbb C)$ be the class with
components $\nu_{\varphi,\alpha}$. Write
$u_a=c_a^{T,\mathrm{top}}(K^{\bullet,\mathrm{an}})$.

The pullback along $r$ of the comparison complex in
\cite[(4.4)]{BrasseletSuwa2000} is $K^{\bullet,\mathrm{an}}$.
Its localized Chern classes are the $u_a$, by the comparison with
Iversen's localized Chern character in
\cite[\S~1(C)]{BrasseletSuwa2000} and the relative Newton identities.
Pulling back \cite[(4.6)]{BrasseletSuwa2000} along $r$ and evaluating
each Chern monomial by the relative cup product therefore gives
\[
\pi^*\rho_\varphi-\nu_\varphi
=
\sum_\lambda A_\lambda\bigl(c(Q^{\mathrm{an}})\bigr)
\smile u_{a_{\lambda,1}}\smile\cdots\smile u_{a_{\lambda,t_\lambda}}
\ \ \text{in} \ \ H^{2d}_T(Y;\mathbb C).
\]

Let $\epsilon\colon H^*_T(Y;\mathbb C)\to H^*(Y;\mathbb C)$ be the
map forgetting the support. For classes $u,v$ supported on $T$,
naturality of the relative cup product gives
$u\smile v=\epsilon(u)\smile v$. Since
$\epsilon(u_a)=c_a(e^{\mathrm{an}})$, we may retain the support on
the first Chern factor in each monomial and forget it on the others.
Thus the $\lambda$-summand is
$P_\lambda^{\mathrm{top}}\smile u_{a_{\lambda,1}}$, where
$P_\lambda^{\mathrm{top}}$ is the cohomology class of $P_\lambda$.
Applying Alexander duality and Corollary~~\ref{cor:localized-chern-comparison}, we obtain
\[
\begin{aligned}
\operatorname{AD}_{Y,T}(\pi^*\rho_\varphi-\nu_\varphi)
&= \sum_\lambda j^*P_\lambda^{\mathrm{top}}
   \cap \operatorname{AD}_{Y,T}(u_{a_{\lambda,1}}) \\
&= \operatorname{cl}_T\!\left(
   \sum_\lambda P_{\lambda,T}\cap
   \Gamma_{a_{\lambda,1}}(K^\bullet)\right)
 = \operatorname{cl}_T(D_\varphi),
\end{aligned}
\]
where the second equality uses the compatibility of the cycle class
map with Chern class operators. Taking the summand supported on
$T_\alpha$ gives
\[
\operatorname{cl}_{T_\alpha}(D_{\varphi,\alpha})
=
\operatorname{AD}_{Y,T_\alpha}
\bigl(\pi^*\rho_{\varphi,\alpha}-\nu_{\varphi,\alpha}\bigr).
\]
\end{proof}

\begin{remark}
\label{rem:independence-localized-factor}
The definition of $D_\varphi$ is independent of which Chern
class factor is localized in each monomial.
Since $Y$ is smooth, its diagonal defines a refined
intersection product
\[
\alpha\star_Y\beta
:=
\Delta_Y^!(\alpha\times\beta)
\in \mathrm{CH}_{n-a-b}(T)_{\mathbb Q},
\]
for $\alpha\in\mathrm{CH}_{n-a}(T)_{\mathbb Q}$, $\beta\in\mathrm{CH}_{n-b}(T)_{\mathbb Q}$ and $T\times_YT\simeq T$.
By compatibility with proper pushforward and
$j_*\Gamma_a=c_a(e)\cap[Y]$, we have
\[
c_b(e)|_T\cap\Gamma_a
=
\Gamma_a\star_Y\Gamma_b
=
c_a(e)|_T\cap\Gamma_b;
\]
see \cite[\S8.1]{Fulton1998}.
Multiplying this identity by the remaining Chern class
operators shows that localizing any other factor
gives the same contribution to
$D_\varphi$ in $\mathrm{CH}_{n-d}(T)_{\mathbb Q}$.
\end{remark}

Let $r_\alpha\colon T_\alpha\to T_\alpha^\nu$ and
$\pi^\nu_\alpha\colon T_\alpha^\nu\to S_\alpha$ be the restrictions
of $r$ and $\pi^\nu$, and put $p_\alpha=\pi^\nu_\alpha\circ r_\alpha$.
Denote by
\[
\operatorname{NRes}_{\varphi,\alpha}
\in H_{2n-2d}(T_\alpha^\nu;\mathbb C)
\]
the image of $\nu^\nu_{\varphi,\alpha}$ under the Alexander
homomorphism for $X^\nu$. This is the sum of the Nash residues
on the connected components of $T_\alpha^\nu$.

\begin{corollary}\label{cor:bb-nash-algebraicity}
Let $\varphi$ be a homogeneous characteristic polynomial with
rational coefficients of degree $q<d\leq n$. For every connected
component $S_\alpha$ of $\operatorname{Sing}(\mathcal F)$, we have
\[
\operatorname{Res}_{\varphi}(\mathcal F,S_\alpha)
-
(\pi^\nu_\alpha)_*\operatorname{NRes}_{\varphi,\alpha}
=
\operatorname{cl}_{S_\alpha}\bigl((p_\alpha)_*D_{\varphi,\alpha}\bigr).
\]
Moreover, the Baum--Bott residue admits a Chow representation on
$S_\alpha$ if and only if the Nash residue admits a Chow
representation on $T_\alpha^\nu$.
\end{corollary}

\begin{proof}
Since $\pi_*[Y]=[X]$ and $r_*[Y]=[X^\nu]$, the projection formula
for cap products with support gives
\[
\begin{aligned}
(p_\alpha)_*\operatorname{AD}_{Y,T_\alpha}
   (\pi^*\rho_{\varphi,\alpha})
&= \operatorname{Res}_{\varphi}(\mathcal F,S_\alpha),\\
(r_\alpha)_*\operatorname{AD}_{Y,T_\alpha}(\nu_{\varphi,\alpha})
&= \operatorname{NRes}_{\varphi,\alpha}.
\end{aligned}
\]
Here the second equality uses
$\nu_{\varphi,\alpha}=r^*\nu^\nu_{\varphi,\alpha}$.
Pushing the second identity of the preceding proposition forward
by $p_\alpha$ and using the compatibility of the cycle class map
with proper pushforward proves the displayed formula. In particular,
if $Z_\alpha$ represents the Nash residue, then
$(\pi^\nu_\alpha)_*Z_\alpha+(p_\alpha)_*D_{\varphi,\alpha}$
represents the Baum--Bott residue.

Conversely, let $R_\alpha\in\mathrm{CH}_{n-d}(S_\alpha)_{\mathbb C}$
represent the Baum--Bott residue. Since $X$ and $Y$ are smooth of
the same dimension, the refined pullback along $\pi$, formed over
$S_\alpha$, gives a class
$\pi^!_{S_\alpha}R_\alpha\in\mathrm{CH}_{n-d}(T_\alpha)_{\mathbb C}$.
Here we identify the Chow group of $Y\times_X S_\alpha$ with that
of its reduction $T_\alpha$. By the compatibility of refined
intersections with the cycle class map
\cite[\S~8.1 and Proposition~19.2]{Fulton1998},
\[
\operatorname{cl}_{T_\alpha}(\pi^!_{S_\alpha}R_\alpha)
=
\operatorname{AD}_{Y,T_\alpha}(\pi^*\rho_{\varphi,\alpha}).
\]
Subtracting the identity of the preceding proposition shows that
$\pi^!_{S_\alpha}R_\alpha-D_{\varphi,\alpha}$ represents
$\operatorname{AD}_{Y,T_\alpha}(\nu_{\varphi,\alpha})$.
Its pushforward
\[
(r_\alpha)_*
\bigl(\pi^!_{S_\alpha}R_\alpha-D_{\varphi,\alpha}\bigr)
\in \mathrm{CH}_{n-d}(T_\alpha^\nu)_{\mathbb C}
\]
therefore represents the Nash residue.
\end{proof}

Let
\(
I_q=(c_{q+1},c_{q+2},\ldots)
\subset\mathbb Q[c_1,c_2,\ldots].
\)
A polynomial belongs to $I_q$ precisely when it vanishes
identically on rank-$q$ vector bundles.
For each $\alpha$, let
\[
\mu_{\alpha,d}:
\CH_{n-d}(S_\alpha)_{\mathbb Q}
\longrightarrow \CH^d(X)_{\mathbb Q}
\]
be proper pushforward along $i_\alpha$, followed by the
inverse of $\cap[X] \colon \CH^d(X)_{\mathbb Q} \to \CH_{n-d}(X)_{\mathbb Q}$.

\begin{theorem}\label{thm:chow-lifts-Iq}
Let $\varphi\in I_q$ be homogeneous of degree $q<d\leq n$. Then the
classes
\[
R_{\varphi,\alpha}:=(p_\alpha)_*D_{\varphi,\alpha}
\in \mathrm{CH}_{n-d}(S_\alpha)_{\mathbb Q}
\]
are Chow representations of the Baum--Bott residues
$\operatorname{Res}_{\varphi}(\mathcal F,S_\alpha)$ and satisfy
\[
\sum_\alpha \mu_{\alpha,d}(R_{\varphi,\alpha})
=
\varphi(N_{\mathcal F})
\ \ \text{in} \ \ \mathrm{CH}^d(X)_{\mathbb Q}.
\]
\end{theorem}

\begin{proof}
Since $\varphi\in I_q$, its restriction to rank-$q$ bundles is
identically zero as an invariant polynomial. Hence the characteristic
and transgression forms defining the localized class of $Q^\nu$
vanish, so $\operatorname{NRes}_{\varphi,\alpha}=0$ for every
$\alpha$. The comparison formula of
Corollary~\ref{cor:bb-nash-algebraicity} therefore gives
\[
\operatorname{Res}_{\varphi}(\mathcal F,S_\alpha)
=
\operatorname{cl}_{S_\alpha}\bigl((p_\alpha)_*D_{\varphi,\alpha}\bigr).
\]

For the global identity, the tautological quotient $Q$ has rank $q$,
so $\varphi(Q)=0$ in $\mathrm{CH}^d(Y)_{\mathbb Q}$. The first
identity of Proposition~3.3 and the projection formula give
\[
\begin{aligned}
\sum_\alpha (i_\alpha)_*R_{\varphi,\alpha}
&= \pi_*j_*D_\varphi \\
&= \pi_*\bigl(\varphi(\xi)\cap[Y]\bigr) \\
&= \varphi(N_{\mathcal F})\cap[X],
\end{aligned}
\]
where the last equality uses
$\varphi(\xi)=\pi^*\varphi(N_{\mathcal F})$ and $\pi_*[Y]=[X]$.
The definition of $\mu_{\alpha,d}$ now gives the asserted equality
in $\mathrm{CH}^d(X)_{\mathbb Q}$.
\end{proof}

\subsection{Codimension-one foliations}\label{subsec:3.3}

Let $\mathcal F$ be a codimension-one foliation on $X$, and denote
\(
\ell=c_1(N_\F)\in H^2(X;\mathbb C).
\)
The preceding subsection gives Chow representations for the residues
of every Chern monomial containing a factor $c_j$ with $j\geq 2$.
In degree $d\geq 2$, the remaining monomial is $c_1^d$. We show
that its residue on each connected component of the singular locus
is determined by the residue of $c_1^2$ on that component.

Write $S=\operatorname{Sing}(\mathcal F)$, let
$i_\alpha\colon S_\alpha\hookrightarrow X$ denote the inclusion
of a connected component, and put
$\rho_{d,\alpha}:=\rho_{c_1^d,\alpha}$.

\begin{proposition}\label{prop:c1-powers-localized}
For every connected component $S_\alpha$ and every $2\leq d\leq n$,
\[
\rho_{d,\alpha}
=
\ell^{d-2}\smile\rho_{2,\alpha}
\qquad\text{in }H^{2d}_{S_\alpha}(X;\mathbb C).
\]
Consequently,
\[
\operatorname{Res}_{c_1^d}(\mathcal F,S_\alpha)
=
(i_\alpha^*\ell)^{d-2}\cap
\operatorname{Res}_{c_1^2}(\mathcal F,S_\alpha).
\]
\end{proposition}

\begin{proof}
We first work with support in $S$. Set $U=X\setminus S$ and
denote by $\rho_d\in H^{2d}_S(X;\mathbb C)$ the localized
Baum--Bott class for $c_1^d$. On $U$ we have
$N_{\mathcal F}|_U\simeq N_\F^{**}|_U$. Choose a basic connection
$\nabla_0$ of type $(1,0)$ on this line bundle extending the
Bott partial connection.

Choose a finite locally free resolution of $N_{\mathcal F}$.
Over $U$, take connections on the resolution compatible with
$\nabla_0$, and over $X$ take arbitrary smooth connections on
its terms. The resulting determinant connections are $\nabla_0$
over $U$ and a connection $\nabla_1$ on $N_\F^{**}$ over $X$.
The first Chern form of the resolution is the first Chern form
of its determinant. This also holds along the path of connections
used to compute the transgression, so both terms of the localized
cocycle for $c_1^d$ can be computed on $N_\F^{**}$;
see \cite[\S~3.5]{CorreaSuwaFunctoriality}.

Write
\[
\Theta_j=\frac{\sqrt{-1}}{2\pi}F_{\nabla_j}
\ \ \textnormal{and} \ \ 
a=\frac{\sqrt{-1}}{2\pi}(\nabla_1|_U-\nabla_0).
\]
In the relative de Rham complex of $(X,U)$, the class $\rho_d$ is
represented by $(\Theta_1^d,T_d)$, where the transgression along
the affine path from $\nabla_0$ to $\nabla_1|_U$ is
\[
T_d
=
d\int_0^1 a\wedge\bigl((1-t)\Theta_0+t\Theta_1\bigr)^{d-1}\,dt
=
a\wedge\sum_{j=0}^{d-1}\Theta_0^j\wedge\Theta_1^{d-1-j}.
\]
Here $\Theta_1$ is restricted to $U$ in the formula for $T_d$.
Bott vanishing for the basic connection gives $\Theta_0^2=0$.
Hence
\[
T_d
=
a\wedge(\Theta_1+\Theta_0)\wedge\Theta_1^{d-2}
=
T_2\wedge\Theta_1^{d-2}.
\]
Thus the relative cocycle for $\rho_d$ is the product of the
cocycle for $\rho_2$ with the global closed form $\Theta_1^{d-2}$.
Since $\Theta_1$ represents $\ell$, we obtain
$\rho_d=\ell^{d-2}\smile\rho_2$. Taking the summand supported
on $S_\alpha$ gives the first identity. The second follows from
the compatibility of Alexander duality with the action of
$H^*(X;\mathbb C)$.
\end{proof}

By Lemma~3.1, the degree-two residue has a unique Chow
representation on each $S_\alpha$. Explicitly,
\[
R_{c_1^2,\alpha}
=
\sum_Z \operatorname{BB}_{c_1^2}(\mathcal F,Z)[Z]
\in \mathrm{CH}_{n-2}(S_\alpha)_{\mathbb C},
\]
where the sum runs over the irreducible components $Z$ of
$S_\alpha$ of codimension two in $X$. If
$\operatorname{codim}_X S_\alpha>2$, the sum is empty and
$H_{2n-4}(S_\alpha;\mathbb C)=0$. In this case
Proposition~\ref{prop:c1-powers-localized} gives
$\operatorname{Res}_{c_1^d}(\mathcal F,S_\alpha)=0$
for every $2\leq d\leq n$.
For $2<d\leq n$, define
\[
R_{c_1^d,\alpha}
:=
c_1(i_\alpha^*L)^{d-2}\cap R_{c_1^2,\alpha}
\in \mathrm{CH}_{n-d}(S_\alpha)_{\mathbb C}.
\]
The cycle class map commutes with the first Chern class operator
of the line bundle $i_\alpha^*N_\F^{**}$. Therefore
Proposition~\ref{prop:c1-powers-localized} gives
\[
\operatorname{cl}_{S_\alpha}(R_{c_1^d,\alpha})
=
(i_\alpha^*\ell)^{d-2}\cap
\operatorname{Res}_{c_1^2}(\mathcal F,S_\alpha)
=
\operatorname{Res}_{c_1^d}(\mathcal F,S_\alpha).
\]

\begin{corollary}\label{cor:codim-one-all-polynomials}
Let $\mathcal F$ be a codimension-one foliation on $X$, and let
$\varphi$ be a homogeneous characteristic polynomial with complex
coefficients of degree $2\leq d\leq n$. Then for every connected
component $S_\alpha$ of $\operatorname{Sing}(\mathcal F)$, the
Baum--Bott residue $\Res_{\varphi}(\F,S_\alpha)$ is algebraic.
\end{corollary}

\begin{proof}
The residue depends linearly on $\varphi$, so it suffices to
consider a Chern monomial. For $c_1^d$, the class
$R_{c_1^d,\alpha}$ constructed above represents the Baum--Bott
residue. Every other monomial of degree $d$ belongs to
$I_1=(c_2,c_3,\ldots)$, so Theorem~\ref{thm:chow-lifts-Iq} 
gives a Chow representation of its Baum--Bott residue.
\end{proof}

\section{The Chow discrepancy}

Theorem~\ref{thm:chow-lifts-Iq} constructs a Chow representation which satisfies an equation equivalent to the Baum--Bott identity \eqref{eq:classical-baum-bott} in the Chow group when the polynomial $\varphi \in I_q = (c_{q+1},c_{q+2},\dots)$. For other characteristic polynomials, it remains to determine whether the global identity holds in the Chow ring. For example the polynomial given by powers of the first Chern class in Proposition~\ref{prop:c1-powers-localized}. We now study this question for codimension-one foliations and  $\varphi = c_1^2$.

For a codimension-one foliation $\F$ on $X$, we denote the pushforward of the Chow representations by
\[
\sum_{\alpha}\mu_*R_{c_1^2,\alpha} = \sum_i \BB_{c_1^2}(\F,Z_i)[Z_i] \ \textnormal{ in } \ \CH^2(X)_\C,
\]
where the right sum is over all irreducible components of $\Sing(\F)$ with codimension two on $X$, and $\BB_{c_1^2}(\F,Z_i) \in \C$ is the Baum--Bott coefficient considered in ~\eqref{eq: BB coef}.
We define the \emph{Chow discrepancy of $\F$} by
\begin{equation}\label{eq: chow-discrepancy}
    \Delta(\F) := c_1(N_\F)^2 - \sum_i\BB_{c_1^2}(\F,Z_i)[Z_i] \in \CH^2(X)_{\C}.
\end{equation}
In particular, we have $\cl_X(\Delta(\F)) = 0$.
For powers of the first Chern class, the representatives constructed
in Subsection~\ref{subsec:3.3} give an explicit discrepancy in every degree.
Let $2< d\leq n$.
On each connected component $S_\alpha$ of the singular locus, take
\[
   R_{c_1^d,\alpha}
   =c_1(i_\alpha^*N_\F)^{d-2}\cap R_{c_1^2,\alpha},
\]
where $R_{c_1^2,\alpha}$ is the unique Chow representative of the
degree-two residue. By the projection formula,
\[
\begin{aligned}
   \Delta_{c_1^d,R}(\mathcal F)
   &:=c_1(N_\F)^d-\sum_\alpha\mu_{\alpha,d}(R_{c_1^d,\alpha})\\
   &=c_1(N_\F)^{d-2}
     \left(c_1(N_\F)^2-\sum_\alpha
       \mu_{\alpha,2}(R_{c_1^2,\alpha})\right)\\
   &=c_1(N_\F)^{d-2}\Delta(\mathcal F)
     \ \ \text{in }\mathrm{CH}^d(X)_{\mathbb C}.
\end{aligned}
\]
Thus, for this choice of representatives, the discrepancy for every
power of $c_1$ is determined by the case $d=2$.

\begin{remark}
Let $\varphi$ be a homogeneous characteristic polynomial of degree
$d$, and suppose that its Baum--Bott residue on each connected
component $S_\alpha$ of $\operatorname{Sing}(\F)$ admits a Chow
representation $R_{\varphi,\alpha}$. For the choice
$R=(R_{\varphi,\alpha})_\alpha$, define the discrepancy by
\[
   \Delta_{\varphi,R}(\F)
   :=\varphi(N_{\F})
     -\sum_\alpha\mu_{\alpha,d}(R_{\varphi,\alpha})
   \in \mathrm{CH}^d(X)_{\mathbb C}.
\]
By the Baum--Bott formula, this class is homologically trivial.
Any other choice $R'$ is of the form
$R'_{\varphi,\alpha}=R_{\varphi,\alpha}+K_\alpha$, with
$K_\alpha\in\ker(\cl_{S_\alpha})$. Hence
\[
   \Delta_{\varphi,R'}(\F)-\Delta_{\varphi,R}(\F)
   =-\sum_\alpha\mu_{\alpha,d}(K_\alpha).
\]
Since the $K_\alpha$ can be chosen independently, the discrepancy
is independent of $R$ if and only if
$\mu_{\alpha,d}(\ker(\cl_{S_\alpha}))=0$ for every $\alpha$.
In particular, for a codimension-one foliation and $\varphi=c_1^2$,
Lemma~\ref{lem: trivial-cases} gives a unique Chow representation on each component,
so we write simply $\Delta(\F)$.
\end{remark}

\begin{proposition}
\label{prop:bb-blowup}
Let $\mathcal F$ be a codimension-one foliation on $X$. Let
\[
\pi\colon X_r\xrightarrow{\pi_r}X_{r-1}
\longrightarrow\cdots
\longrightarrow X_1\xrightarrow{\pi_1}X_0=X,
\]
be a finite composition of blow-ups along smooth closed centers
of positive codimension. Set $\mathcal F_0=\mathcal F$, and let
$\mathcal F_i$ be the pullback of
$\mathcal F_{i-1}$ under $\pi_i$. Then
\[
\pi_*\Delta(\mathcal F_r)=\Delta(\mathcal F).
\]
\end{proposition}

\begin{proof}
First suppose that
\(
\pi\colon Y=\operatorname{Bl}_W(X)\longrightarrow X
\)
is the blow-up of a smooth irreducible center $W$, and write
$\mathcal G$ for the saturated transform of $\mathcal F$.
If $\operatorname{codim}_X W=1$, then $W$ is an effective
Cartier divisor and $\pi$ is an isomorphism. Hence we may
assume that $\operatorname{codim}_X W\geq 2$.

The foliations and their normal sheaves agree outside the
center $W$. Hence, by localization for Chow groups,
\[
\Gamma:=\pi_*\Delta(\mathcal G)-\Delta(\mathcal F)
\in
\operatorname{im}\bigl(
\mathrm{CH}_{n-2}(W)_{\mathbb C}
\longrightarrow \mathrm{CH}_{n-2}(X)_{\mathbb C}
\bigr).
\]
If $\operatorname{codim}_XW>2$, localization gives $\Gamma=0$.
Otherwise, $\Gamma=a[W]$, and pairing with an ample class
to the power $n-2$ gives $a=0$. Therefore $\pi_*\Delta(\mathcal{G}) = \Delta(\F)$.

A smooth closed center has finitely many pairwise disjoint
smooth irreducible components. Blowing up their union is
equivalent to blowing up the components successively.
Applying the single blow-up identity at each stage gives
\(
(\pi_i)_*\Delta(\mathcal F_i)
=
\Delta(\mathcal F_{i-1}).
\)
Functoriality of proper pushforward gives
\[
\pi_*\Delta(\mathcal F_r)
=
\Delta(\mathcal F).
\]
\end{proof}

The above result allows us to compute the discrepancy on a suitable birational model. We will use such a model to study the Chow discrepancy for foliations defined by closed rational one-forms. 

\subsection{Closed rational one-forms and the discrepancy formula}\label{subsec:closed-one-form}

We will focus on codimension-one foliations defined by closed rational one-forms. Let $\eta$ be a closed rational one form on $X$, and let $D$ be a divisor such that
\[
D = -\mathrm{div}(\eta) = (\eta)_\infty - (\eta)_0.
\]
The tangent sheaf of the codimension-one foliation defined by $\eta$ is given by
\[
T_{\mathcal F_\eta} = \ker\!\left(T_X \xrightarrow{\,v\mapsto\eta(v)\,}\mathcal O_X (D)\right).
\]
Let $D_i$ be an irreducible divisor with $x \in D_i$ a smooth point. Let $p\colon \overline{\mathbb{D}} \hookrightarrow X$ be an embedding of a closed disc, holomorphic in the interior and intersecting $D_i$ transverselly at $x$. We define the residue of $\eta$ along $D_i$ by
\[
\Res_{D_i}(\eta) = \frac{1}{2\pi i} \int_{\partial \mathbb{D}} p^* \eta \in \C.
\]
Closedness of $\eta$ and irreducibility of the divisor make the residue well-defined. Let $D_i$ be the irreducible components of $(\eta)_{\infty}$, with pole order $p_i \in \mathbb{Z}_{>0}$.
We define the residue divisor of $\eta$ by
\[
    \operatorname{Res}(\eta)
    := \sum_i \operatorname{Res}_{D_i}(\eta)\,D_i
    \in \operatorname{Div}(X)\otimes_{\mathbb Z}\mathbb C.
\]
We use the same notation for its class in
$\mathrm{CH}^1(X)_{\mathbb C}$. By Weil's theorem \cite[Theorem 9.3.1]{PereiraClosedForms} this class is zero in cohomology, therefore $\Res(\eta) \in \Pic^0(X)_{\C}$. Let $\lambda_i := \Res_{D_i}(\eta)$. We define the divisor modulo rational equivalence
\[
B_{\eta} : = \sum_{p_i = 1} \lambda_i^{-1} [D_i],
\]
where the sum is over all irreducible components of $(\eta)_{\infty}$ with simple poles, and therefore $\lambda_j \neq 0$. If $\eta$ has no simple poles, then $B_{\eta} = 0$.

\begin{theorem}
\label{thm:closed-rational-discrepancy}
Let $X$ be a smooth projective variety over $\C$ of dimension $n\geq 2$, and let $\eta$ be a nonzero closed rational one-form on $X$.
Then
\[
\Delta(\mathcal{F}_\eta)
=
\Res(\eta) B_\eta.
\]
In particular,
\(
\Delta(\mathcal{F}_\eta)
\in
\operatorname{Pic}^0(X)_{\mathbb{C}}
\cdot\operatorname{CH}^1(X)_{\mathbb{C}} \subset \CH^2(X)_{\alg,\C}.
\)
\end{theorem}

We prove Theorem~\ref{thm:closed-rational-discrepancy} in Subsection~\ref{subsec:residue-calculations-proof}. Example~\ref{ex:elliptic} shows that the discrepancy can be nonzero modulo rational equivalence.

\begin{remark}
Let $D$ be a reduced divisor on $X$. A rational one-form $\eta$ is
logarithmic along $D$ if it is holomorphic on $X\setminus D$ and
both $\eta$ and $d\eta$ have at most simple poles along $D$.
Suppose that $\eta$ is closed. Near a point of $D$,
let $f_1,\ldots,f_r$ be reduced local equations for its irreducible
analytic branches, and let $\lambda_i$ be the residue of $\eta$
along the branch $(f_i=0)$. If $\eta$ is logarithmic along $D$,
then locally
\[
   \eta=\sum_{i=1}^{r}\lambda_i\frac{df_i}{f_i}+dh
\]
for some holomorphic function $h$. Indeed, subtracting the logarithmic
terms removes the pole at a general smooth point of each branch.
The difference is therefore holomorphic, and it is locally exact
by closedness.
When $D$ has simple normal crossings, every global logarithmic
one-form along $D$ is closed since $X$ is smooth and projective.
In this case the $f_i$ may be chosen as part of a local coordinate
system.
\end{remark}

We recall Pereira’s definition of simple singularities for a closed rational one-form \cite[Definition 9.4.4]{PereiraClosedForms}. Let $\eta$ be a closed rational one-form on $X$. Locally, on holomorphic coordinates $x_i$, one may write the form as
\[
\eta = \sum_i \lambda_i\frac{dh_i}{h_i} + d \left( \frac{f}{g} \right)
\]
where the $h_i$ are reduced local equations of the branches with nonzero residue and f,g are relatively prime.
We define \emph{base locus} of $\eta$ as the analytic subspace of $X$ defined by the ideal $\mathfrak{B}(\eta) \subset \mathcal{O}_X$ locally generated by $f,g$, see \cite[Subsection 9.4.2]{PereiraClosedForms}. The \emph{irregular divisor} $\textnormal{Irr}(\eta)$ is defined by the difference between the polar divisor $(\eta)_{\infty}$ and its reduced support. We say that $\eta$ is resonant at a point $p$ if there exists irreducible components $D_1,\dots,D_k$ on the support of $\Res(\eta)$ containing $p$ none of them contained in the support of $\textnormal{Irr}(\eta)$ and positive integers $n_1,\dots,n_k$ such that
\[
\sum_{i=1}^k n_i \Res_{D_i}(\eta) = 0.
\]
We say that $\eta$ is \emph{non-resonant}, if it is not resonant at every point $p\in X$, see \cite[Subsection 9.4.4]{PereiraClosedForms}.

\begin{definition}
    Let $\eta$ be a closed rational one-form on $X$. We say that $\eta$ has \emph{simple singularities} if $\mathfrak B(\eta) = \mathcal{O}_X$, the form is non-resonant and the support of $(\eta)_{\infty}$ is a simple normal crossing divisor.
\end{definition}

\begin{lemma}[\cite{PereiraClosedForms}, Proposition 9.4.5]\label{lem:simple-model}
Let $\eta$ be a
nonzero closed rational one-form on $X$. There exists a finite
composition of blow-ups along smooth closed algebraic centers
\[
    \pi \colon Y \longrightarrow X
\]
such that $\omega =\pi^*\eta$ has simple singularities.
\end{lemma}
\begin{proof}
    It follows from the proof of \cite[Proposition 9.4.5]{PereiraClosedForms}. First take principalization, see \cite[Theorem 1.0.1]{WlodarczykResolution}, of the ideal $\mathfrak{B}(\eta)$ giving a sequence of blow-ups. Then take an embedded resolution, see \cite[Theorem 1.0.2]{WlodarczykResolution}, of the support of the polar divisor of the closed one-form on the first sequence of blowups. Finally, we eliminate the resonances using \cite[Theorem 2]{FernandezDuque2015}.
\end{proof}

\begin{corollary}\label{cor:zero-locus-irregular-divisor}
Let $\omega=\pi^*\eta$ be as in Lemma~\ref{lem:simple-model}.
Set $D$ the reduced support of $(\omega)_\infty$ and
$U=Y\setminus D$. Then the following hold.
\begin{enumerate}
    \item The reduced zero locus $K=Z(\omega|_U)_{\mathrm{red}}$ is
    closed in $Y$ and hence projective. For each connected component
    $K_\alpha$ of $K^{\mathrm{an}}$, there exist an analytic neighborhood
    $V_\alpha\subset U^{\mathrm{an}}$ of $K_\alpha$ and a holomorphic
    function $f_\alpha$ on $V_\alpha$ such that
    \[ 
        df_\alpha=\omega|_{V_\alpha}
    \ \ \textnormal{and} \ \ f_\alpha|_{K_\alpha}=0.
    \]
    The neighborhoods $V_\alpha$ may be chosen pairwise disjoint.
    \item The irregular divisor $I =\textnormal{Irr}(\omega)$ satisfies
    \(
        \mathcal O_Y(I)|_I \simeq \mathcal O_I .
    \)
\end{enumerate}
\end{corollary}

\begin{proof}
By \cite[Corollary~9.4.7]{PereiraClosedForms}, the analytic closure of
$K^{\mathrm{an}}$ in $Y^{\mathrm{an}}$ is disjoint from $D^{\mathrm{an}}$.
Since $K$ is closed in $U$, it follows that $K^{\mathrm{an}}$ is a
closed analytic subset of $Y^{\mathrm{an}}$. Hence $K$ is a closed
subvariety of $Y$ by Chow's theorem.
By \cite[Lemma 9.4.8]{PereiraClosedForms}, $\omega$ admits a holomorphic primitive $f_\alpha$ on a neighborhood of each connected component $K_\alpha$. Since $df_\alpha$ vanishes along $K_\alpha$, the restriction of $f_\alpha$ to $K_\alpha$ is constant.

By \cite[Proposition~9.4.9]{PereiraClosedForms}, the line bundle
$\mathcal O_Y(I)|_I$ is analytically trivial. Since $I$ is a
projective scheme, GAGA gives the asserted algebraic trivialization.
\end{proof}

Since $df_\alpha=\omega|_{V_\alpha}$, the reduced zero fiber of
$f_\alpha$ is smooth outside $K_\alpha$. By embedded resolution \cite[Theorem A]{BierstoneMilman}
near each compact $K_\alpha$, after finitely many further blow-ups
with smooth centers lying over $K$, we may assume that
$\operatorname{div}(f_\alpha)$ has simple normal crossings support
for every $\alpha$. Each center is a compact analytic subvariety
of the corresponding projective model and is therefore algebraic.

We keep the notation $Y$, $\omega$, and $f_\alpha$ on this
model, writing $V_\alpha$ for the inverse image of the original
neighborhood and $K_\alpha$ for the reduced inverse image of the
original component. The modification is an isomorphism away from
$K$, and the pulled-back form vanishes over $K$. Hence, on each
$V_\alpha$,
\[
    df_\alpha=\omega|_{V_\alpha}
    \ \ \textnormal{and} \ \
    \operatorname{Crit}(f_\alpha)_{\mathrm{red}}
    = K_\alpha \subseteq f_\alpha^{-1}(0).
\]

\subsection{Residue calculations and proof of the formula}\label{subsec:residue-calculations-proof}
By Proposition~\ref{prop:bb-blowup}, we may compute the Chow discrepancy on a
birational model obtained by a sequence of blow-ups along smooth
centers and then push it forward to $X$. We first choose a model
on which the pullback of $\eta$ has simple singularities. On this model,
the zero locus of the form away from its poles is disjoint from
the polar divisor even after taking its closure. This allows us
to compute the two contributions to the Baum--Bott cycle separately.

\begin{lemma}
\label{lem:closed-form-zero-cancellation} 
Define
\[
\mathcal B_{\mathrm{zero}}
:=
\sum_{
Z\not\subset\operatorname{Supp}(\omega)_{\infty}}
\operatorname{BB}_{c_1^2}(\mathcal F_\omega,Z)[Z]
\in \mathrm{CH}^2(Y)_{\mathbb C},
\]
where the sum is over all irreducible components of $\Sing(\F_\omega)$ with codimension two on $Y$.
Then 
\[
\mathcal B_{\mathrm{zero}}=[(\omega)_0]^2.
\]
\end{lemma}

\begin{proof}
On a neighborhood $V_\alpha$ provided by Corollary ~\ref{cor:zero-locus-irregular-divisor}, write
\[
    F_\alpha:=\operatorname{div}(f_\alpha)
       =\sum_a m_aE_a,
    \ \ \textnormal{with} \ \ m_a\geq1.
\]
We may assume that
only finitely many components occur. Since the support of $F_\alpha$
is simple normal crossing, locally, after absorbing a holomorphic
unit into one coordinate, we have
\(
    f_\alpha=x_1^{m_1}\cdots x_r^{m_r}.
\)
Consequently, the divisorial zero part of $df_\alpha$ is
\[
    R_\alpha=\sum_{m_a>1}(m_a-1)E_a.
\]

If $m_a>1$, then $E_a$ is contained in the critical set of
$f_\alpha$. It is therefore a compact closed analytic hypersurface
of $Y^{\mathrm{an}}$, and hence a smooth algebraic divisor on $Y$. Likewise, for $a\ne b$, the intersection
$E_a\cap E_b$ is contained in the critical set and is a compact
closed analytic subset of $Y^{\mathrm{an}}$.
It thus defines an algebraic codimension-two cycle
\[
    \Gamma_{ab}:=[E_a\cap E_b].
\]
The intersections are transverse, so their irreducible components
occur with multiplicity one. We do not require that the
components with $m_a=1$ extend to algebraic divisors on $Y$.

Fix $a$ with $m_a>1$. Since $F_\alpha$ is principal on $V_\alpha$,
restricting its associated line bundle to $E_a$ gives an analytic
trivialization of
\[
    N_{E_a/Y}^{\otimes m_a}
       \otimes
    \mathcal O_{E_a}
       \left(\sum_{b\ne a}m_b(E_a\cap E_b)\right).
\]
Each intersection in this expression is an algebraic Cartier
divisor on the smooth projective variety $E_a$. By \cite{SerreGAGA}, the
trivialization is algebraic. Taking first Chern classes on $E_a$
and pushing forward to $Y$ yields
\[
    m_a[E_a]^2+\sum_{b\ne a}m_b\Gamma_{ab}=0
    \ \ \text{in }\mathrm{CH}^2(Y).
\]

Expanding $[R_\alpha]^2$ and eliminating its self-intersections using
these relations, the coefficient of $\Gamma_{ab}$ becomes
\[
    2(m_a-1)(m_b-1)
    -\frac{(m_a-1)^2m_b}{m_a}
    -\frac{(m_b-1)^2m_a}{m_b}
    =-\frac{(m_a-m_b)^2}{m_am_b}.
\]
When either multiplicity equals one, the terms with the corresponding
zero factor are absent. Thus no self-intersection class is assigned
to a component with multiplicity one, and
\[
    [R_\alpha]^2
      =-\sum_{a<b}
          \frac{(m_a-m_b)^2}{m_am_b}\,\Gamma_{ab}
    \ \ \text{in }\mathrm{CH}^2(Y)_{\mathbb C}.
\]

At a general point of $E_a\cap E_b$, the transverse first integral
is $x^{m_a}y^{m_b}$. Its saturated differential is
$m_a y\,dx+m_b x\,dy$, whose tangent vector field
$m_bx\,\partial_x-m_ay\,\partial_y$ has eigenvalues $m_b,-m_a$.
The transverse Baum--Bott residue formula therefore gives
\[
    \operatorname{BB}_{c_1^2}
       (\mathcal F_\omega,Z_{ab,\nu})
       =-\frac{(m_a-m_b)^2}{m_am_b}
\]
for every irreducible component $Z_{ab,\nu}$ of $E_a\cap E_b$;
see \cite[Theorem~1.2]{CorreaLourenco2019}. This coefficient is independent of $\nu$. Since $\Gamma_{ab} = \sum_{\nu} [Z_{ab,\nu}]$, we have
\[
\begin{aligned}
\sum_{Z\subset V_\alpha}
\operatorname{BB}_{c_1^2}(\mathcal F_\omega,Z)[Z]
&=
-\sum_{a<b}
\frac{(m_a-m_b)^2}{m_am_b}
\sum_{\nu}[Z_{ab,\nu}]\\
&=
-\sum_{a<b}
\frac{(m_a-m_b)^2}{m_am_b}\,\Gamma_{ab}.
\end{aligned}
\]
At a point belonging to only one component of $F_\alpha$,
the foliation is regular. Intersections of at least three
components have codimension at least three, and there are no critical
points outside the zero fiber in $V_\alpha$.
Thus the total codimension-two Baum–Bott contribution on $V_\alpha$ is $[R_\alpha]^2$. If $m_a=m_b=1$, the transverse first integral is $xy$.
Its differential $y\,dx+x\,dy$ has an isolated singularity
with Baum--Bott residue zero. Thus these intersections are
singular components, but have coefficient zero in
$\mathcal B_{\mathrm{zero}}$ and in the displayed expression for
$[R_\alpha]^2$.

Finally, the neighborhoods $V_\alpha$ are pairwise disjoint, so the
algebraic divisors $R_\alpha$ have disjoint supports. Since
$(\omega)_0=\sum_\alpha R_\alpha$, their mixed intersection products in the Chow ring vanish,
and summing gives
\[
  \mathcal  B_{\mathrm{zero}}
       =\sum_\alpha [R_\alpha]^2
       =[(\omega)_0]^2.
\]
\end{proof}

\begin{lemma}
\label{lem:closed-form-polar-contribution}
Define
\[
\mathcal B_{\mathrm{pole}}
:=
\sum_{Z\subset\operatorname{Supp}(\omega)_{\infty}}
\operatorname{BB}_{c_1^2}(\mathcal F_\omega,Z)[Z]
\in \mathrm{CH}^2(Y)_{\mathbb C},
\]
where the sum is over all irreducible components of $\Sing(\F_\omega)$ with codimension two on $Y$.
Then
\[
[(\omega)_{\infty}]^2-\mathcal B_{\mathrm{pole}}=\Res(\omega)B_{\omega}.
\]
\end{lemma}

\begin{proof}
Write $(\omega)_{\infty} = \sum_i (k_i+1) D_i$ with $k_i \geq0$ and put $\lambda_i = \Res_{D_i}(\omega)$. At a point belonging to only one polar component, the normal form
of Corollary~\ref{cor:zero-locus-irregular-divisor} defines a regular foliation.
Thus every codimension-two polar singular component is an
irreducible component of some $D_i\cap D_j$.
Since the polar support is simple normal crossing, the cycle
$[D_i][D_j]$ is the sum of these components with multiplicity one.
At a general point of $D_i\cap D_j$, the transverse normal form is
\[
    \omega
      =\lambda_i\frac{dx}{x}
       +\lambda_j\frac{dy}{y}
       +d(x^{-k_i}y^{-k_j}).
\]
Write $\beta_{ij}$ for the corresponding Baum--Bott coefficient.
We use the transverse residue formula of
\cite[Theorem~1.2]{CorreaLourenco2019}.

If $k_i=k_j=0$, the saturated tangent vector field is
$\lambda_j x\,\partial_x-\lambda_i y\,\partial_y$.
If $k_i,k_j>0$, putting $t=x^{k_i}y^{k_j}$, it is
\[
    x(\lambda_j t-k_j)\partial_x+y(k_i-\lambda_i t)\partial_y.
\]
Their linearizations have eigenvalues $\lambda_j,-\lambda_i$ and $-k_j,k_i$,
respectively. The nondegenerate residue formula therefore gives
\[
    \beta_{ij}
      =
      \begin{cases}
        \displaystyle
        -\frac{(\lambda_i-\lambda_j)^2}{\lambda_i\lambda_j},
          & k_i=k_j=0,\\[6pt]
        \displaystyle
        -\frac{(k_i-k_j)^2}{k_ik_j},
          & k_i,k_j>0.
      \end{cases}
\]
In the first case, $\lambda_i$ and $\lambda_j$ are nonzero because the
corresponding components are actual simple poles.

Suppose now that $k_i>0$ and $k_j=0$. Then $\lambda_j\ne0$, and a saturated
tangent vector field is
\[
    v=\lambda_j x^{k_i+1}\partial_x+(k_i-\lambda_i x^{k_i})y\partial_y.
\]
Its divergence is
\[
    \operatorname{div}(v)
       =k_i+\bigl((k_i+1)\lambda_j-\lambda_i\bigr)x^{k_i}.
\]
Since $k_i-\lambda_i x^{k_i}$ is a unit near the origin, taking the
Grothendieck residue first in the $y$-variable gives
\[
    \beta_{ij}
       =\operatorname{Res}_{x=0}
          \frac{
            \left[k_i+\bigl((k_i+1)\lambda_j-\lambda_i\bigr)x^{k_i}\right]^2
          }{
            \lambda_j x^{k_i+1}(k_i-\lambda_i x^{k_i})
          }\,dx.
\]
Equivalently, this is the coefficient of $x^{k_i}$ in
\[
    \frac{
      \left[k_i+\bigl((k_i+1)\lambda_j-\lambda_i\bigr)x^{k_i}\right]^2
    }{
      \lambda_j(k_i-\lambda_i x^{k_i})
    }.
\]
Using
\[
    \frac{1}{k_i-\lambda_i x^{k_i}}
       =\frac{1}{k_i}+\frac{\lambda_i}{k_i^2}x^{k_i}+O(x^{2k_i}),
\]
we obtain
\[
    \beta_{ij}
       =\frac{
          2k_i\bigl((k_i+1)\lambda_j-\lambda_i\bigr)+k_i\lambda_i
        }{\lambda_j k_i}
       =2(k_i+1)-\frac{\lambda_i}{\lambda_j}.
\]
The case $k_i=0$ and $k_j>0$ follows by exchanging the indices.
We now compare these coefficients with the intersection products.
Set $H=\{i:k_i>0\}$ and $L=\{i:k_i=0\}$
and write
\[
    P_H=\sum_{i\in H}(k_i+1)[D_i],
    \qquad
    P_L=\sum_{i\in L}[D_i],
    \qquad [(\omega)_{\infty}]=P_H+P_L.
\]
Also put
\[
    A_H=\sum_{i\in H}\lambda_i[D_i],
    \qquad
    A_L=\sum_{i\in L}\lambda_i[D_i],
    \qquad \Res(\omega)=A_H+A_L.
\]

For $i\in H$ and $I = \textnormal{Irr}(\omega)$, the algebraic trivialization
$\mathcal O_Y(I)|_I\simeq\mathcal O_I$ restricts along
$D_i\hookrightarrow I$ to
\(
    \mathcal O_Y(I)|_{D_i}\simeq\mathcal O_{D_i}.
\)
Taking first Chern classes on $D_i$ and pushing forward to $Y$
therefore gives
\(
    [I][D_i]=0,
\)
or explicitly,
\[
    k_i[D_i]^2+
       \sum_{\substack{j\in H\\j\ne i}}
          k_j[D_i][D_j]=0
    \ \ \text{in }\mathrm{CH}^2(Y).
\]

Decompose the polar Baum--Bott cycle as
\[
    \mathcal B_{\mathrm{pole}}=B_{HH}+B_{HL}+B_{LL},
\]
according to whether the two intersecting components belong to
$H$, to different sets, or to $L$.
Eliminating the self-intersections in $P_H^2$ with the preceding
relations, the coefficient of $[D_i][D_j]$, for $i,j\in H$, becomes
\[
    2(k_i+1)(k_j+1)
      -\frac{(k_i+1)^2k_j}{k_i}
      -\frac{(k_j+1)^2k_i}{k_j}
      =-\frac{(k_i-k_j)^2}{k_ik_j}.
\]
These are exactly the coefficients of $B_{HH}$, so
\[
    P_H^2=B_{HH}.
\]

For two simple-pole components, the residue formula gives
\[
\begin{aligned}
    P_L^2-B_{LL}
       &=
       \sum_{i\in L}[D_i]^2
       +\sum_{\substack{i<j\\i,j\in L}}
          \left(
             \frac{\lambda_i}{\lambda_j}
             +\frac{\lambda_j}{\lambda_i}
          \right)[D_i][D_j]\\
       &=
       \left(\sum_{i\in L}\lambda_i[D_i]\right)
       \left(\sum_{j\in L}\lambda_j^{-1}[D_j]\right)
       =A_LB_{\omega}.
\end{aligned}
\]
For a mixed pair $i\in H$, $j\in L$, the coefficient of
$[D_i][D_j]$ in $2P_HP_L-B_{HL}$ is
\[
    2(k_i+1)
       -\left(2(k_i+1)-\frac{\lambda_i}{\lambda_j}\right)
       =\frac{\lambda_i}{\lambda_j}.
\]
Consequently,
\[
    2P_HP_L-B_{HL}
       =\sum_{\substack{i\in H\\j\in L}}
          \frac{\lambda_i}{\lambda_j}[D_i][D_j]
       =A_HB_{\omega}.
\]
Adding the three identities yields
\[
\begin{aligned}
    [(\omega)_{\infty}]^2- \mathcal B_{\mathrm{pole}}
       &=(P_H^2-B_{HH})
         +(2P_HP_L-B_{HL})
         +(P_L^2-B_{LL})\\
       &=A_HB_{\omega}+A_LB_{\omega}
        =\Res(\omega)B_{\omega}.
\end{aligned}
\]
\end{proof}

\begin{proof}[Proof of Theorem \ref{thm:closed-rational-discrepancy}]
Choose the birational model constructed in Subsection~\ref{subsec:closed-one-form}, denoted  $\pi\colon Y\to X$ and put $\omega=\pi^*\eta$, with the notation fixed.
Then we have that
\(
\det(N_{\mathcal{F}_\omega})
\simeq
\mathcal{O}_Y((\omega)_{\infty}-(\omega)_0).
\)
Since $\Supp (\omega)_{\infty} \cap K = \emptyset$ and $\Supp (\omega)_0 \subset K$, the supports of $(\omega)_{\infty}$ and $(\omega)_0$ are disjoint.
Consequently $[(\omega)_{\infty}] [(\omega)_0]=0$ in the Chow ring and
\[
c_1(N_{\mathcal{F}_\omega})^2= c_1(\det N_{\F_\omega})^2 =[(\omega)_{\infty}]^2+[(\omega)_0]^2.
\]
Every codimension-two singular component lies either in the
polar locus or in the zero locus outside it.
Lemmas~\ref{lem:closed-form-zero-cancellation}
and~\ref{lem:closed-form-polar-contribution} therefore give
\[
\begin{aligned}
\Delta(\mathcal{F}_\omega)
&=
[(\omega)_{\infty}]^2+[(\omega)_0]^2
-\mathcal{B}_{\mathrm{pole}}
-\mathcal{B}_{\mathrm{zero}}\\
&=
[(\omega)_{\infty}]^2-\mathcal{B}_{\mathrm{pole}}\\
&=
\Res(\omega) B_\omega.
\end{aligned}
\]
We descend this identity to $X$.
By \cite[Lemma 9.4.3]{PereiraClosedForms} we have
\[
\Res(\omega)=\pi^*\Res(\eta).
\]
Moreover, $\pi$ is an isomorphism over the generic point of each
prime divisor $D\subset X$. Thus $D$ is a simple pole of $\eta$
if and only if its strict transform $\widetilde D$ is a simple pole
of $\omega$, in which case their residues agree. Since
$\pi_*[\widetilde D]=[D]$ and every exceptional prime divisor
has image of codimension at least two in $X$, the definition of
$B_\omega$ gives
\[
    \pi_*B_\omega=B_\eta.
\]
Using Proposition \ref{prop:bb-blowup}
and the projection formula, we obtain
\[
\begin{aligned}
\Delta(\mathcal{F}_\eta)
&=
\pi_*\Delta(\mathcal{F}_\omega)
=
\pi_*\bigl(\pi^*\Res(\eta)\cdot B_\omega\bigr)
=
\Res(\eta)\cdot\pi_*B_\omega\\
&=
\Res(\eta) B_\eta.
\end{aligned} 
\]\qedhere
\end{proof}

\subsection{Consequences and examples}
We first give an example in which the discrepancy is nonzero in the Chow group.

\begin{example}\label{ex:elliptic}
Let $E$ be an elliptic curve and choose $p,q\in E$ with
$p-q$ non-torsion in $\Pic^0(E)$.
By Riemann--Roch there exists a meromorphic differential $\alpha$
on $E$ with precisely two simple poles and residues $1,-1$ at
$p,q$. On $\PP^1$ take
\[
\beta=\frac{dt}{t}+\frac{dt}{t-1},
\]
whose residues at $0,1,\infty$ are $1,1,-2$.
On $S=E\times\PP^1$ define
\[
\eta=\operatorname{pr}_E^*\alpha+
     \operatorname{pr}_{\PP^1}^*\beta.
\]
Write $F_x=\{x\}\times\PP^1$ and $H_t=E\times\{t\}$.
The polar divisor is
\[
(\eta)_{\infty}=F_p+F_q+H_0+H_1+H_\infty.
\]
It has simple normal crossings.
All $H_t$ have the same Chow class $H$, so
\[
\Res(\eta)=F_p-F_q
\ \ \ \textnormal{and} \ \ \
B_\eta=F_p-F_q+\frac32H.
\]
Since $(F_p-F_q)^2=0$, Theorem~\ref{thm:closed-rational-discrepancy} gives
\[
\Delta(\F_\eta)
=\frac32(F_p-F_q)H
=\frac32\bigl([(p,t_0)]-[(q,t_0)]\bigr)\neq0
\]
for any $t_0\in\PP^1$.
Nonvanishing follows by pushing forward to the elliptic curve $E$.
\end{example}

\begin{example}\label{ex:products}
For $n\geq2$, take
\(
X_n=E\times\PP^1\times\PP^{n-2}
\)
and pull back the form of Example~\ref{ex:elliptic} under the
projection $r\colon X_n\to E\times\PP^1$.
The form is logarithmic and its discrepancy is
$r^*\Delta(\F_\eta)$. This is nonzero by the projective bundle
formula.
\end{example}

Theorem~\ref{thm:closed-rational-discrepancy} shows that the discrepancies associated with closed
rational one-forms belong to the product space
$\operatorname{Pic}^0(X)_{\mathbb C}\cdot\operatorname{CH}^1(X)_{\mathbb C}$.
We now show that discrepancies of closed logarithmic one-forms span this space.  

Fix a nontrivial element $L \in \Pic^0(X)$, and let $M$ be a line bundle such that $M$ and $M\otimes L$ are very ample. Denote \(\ell=c_1(L)\) and \(h = c_1(M)\).

\begin{lemma}\label{lem:three-pole-construction}
For general divisors $D_1\in |M\otimes L|$ and $D_2,D_3\in |M|$,
there exists a closed rational one-form $\eta_{L,M}$ on $X$ with
polar divisor $D_1+D_2+D_3$ and residues $1,1,-2$ along
$D_1,D_2,D_3$, respectively. Its Chow discrepancy is
\[
    \Delta(\mathcal F_{\eta_{L,M}})
    = \ell^2+\frac{3}{2}\ell h .
\]
\end{lemma}

\begin{proof}
By Bertini's theorem, the divisors $D_1,D_2,D_3$ are smooth,
irreducible, and distinct, and their union has simple normal
crossings. Set $R=D_1+D_2-2D_3$. Since $\mathcal O_X(R)\simeq L$
and $L\in\operatorname{Pic}^0(X)$, the class of $R$ vanishes in
$H^2(X,\mathbb C)$. By \cite[Proposition~9.5.5]{PereiraClosedForms}, there
exists a closed logarithmic one-form $\eta_{L,M}$ with residue
divisor $R$. Since $X$ is projective, this form is rational.
Its poles are simple and occur precisely along $D_1,D_2,D_3$.
In $\mathrm{CH}^1(X)_{\mathbb C}$ we have
$[D_1]=h+\ell$ and $[D_2]=[D_3]=h$. Hence
\[
    \operatorname{Res}(\eta_{L,M})=\ell
    \ \ \textnormal{and} \ \
    B_{\eta_{L,M}}
    =[D_1]+[D_2]-\tfrac12[D_3]
    =\ell+\tfrac32 h .
\]
The formula follows from Theorem~\ref{thm:closed-rational-discrepancy}.
\end{proof}

Now consider the subspace of $\CH^2(X)_{\C}$ defined by
\[
V_{\mathrm{cl}}(X)
:=
\operatorname{Span}_{\mathbb{C}}
\left\{
\Delta(\mathcal{F}_\eta):
\eta\text{ is a nonzero closed rational one-form}
\right\}.
\]

\begin{corollary}
\label{cor:discrepancy-span}
Let $X$ be a smooth projective variety over
$\mathbb C$ of dimension $n \geq 2$. Then
\[
V_{\mathrm{cl}}(X)
=
\operatorname{Pic}^0(X)_{\mathbb{C}}
\cdot\operatorname{CH}^1(X)_{\mathbb{C}}.
\]
\end{corollary}

\begin{proof}
Theorem~\ref{thm:closed-rational-discrepancy} gives
\(
V_{\mathrm{cl}}(X)
\subseteq
\operatorname{Pic}^0(X)_{\mathbb{C}}
\cdot\operatorname{CH}^1(X)_{\mathbb{C}}.
\)
For the reverse inclusion, fix
$L\in\operatorname{Pic}^0(X)$ and
$B\in\operatorname{Pic}(X)$. Put $\ell = c_1(L)$ and $b = c_1(B)$.
Choose $M$ sufficiently ample that $M$, $ M\otimes L$, $ M \otimes B$, $M \otimes B \otimes L$
are all very ample and write $h=c_1(M)$.
Apply Lemma~\ref{lem:three-pole-construction} using $M$
and $M\otimes B$. It gives
\[
\Delta(\mathcal{F}_{\eta_{L,M}})
=\ell^2+\frac32\ell h \ \ \textnormal{and} \ \
\Delta(\mathcal{F}_{\eta_{L,M\otimes B}})
=\ell^2+\frac32\ell(h+b).
\]
Subtracting yields
\[
c_1(L)c_1(B)
=
\frac23
\left(
\Delta(\mathcal{F}_{\eta_{L,M\otimes B}})
-
\Delta(\mathcal{F}_{\eta_{L,M}})
\right).
\]
Hence every such product belongs to $V_{\mathrm{cl}}(X)$.
\end{proof}

We now characterize the varieties $X$ for which every discrepancy associated
with a closed rational one-form vanishes. By Corollary~\ref{cor:discrepancy-span}, this amounts
to determining when
$\operatorname{Pic}^0(X)_{\mathbb C}\cdot\operatorname{CH}^1(X)_{\mathbb C}$
is zero. 
Let
\[
a_X:\mathrm{CH}_0(X)^0\longrightarrow\operatorname{Alb}(X)(\mathbb C)
\]
be the Albanese homomorphism, where $\mathrm{CH}_0(X)^0$
denotes the subgroup of zero-cycles of degree zero.
We use the following standard consequence of
Picard--Albanese duality and polarization theory.

\begin{lemma} \label{lem:ample-albanese-isogeny}
Let $h=c_1(M)$ for a very ample
line bundle $M$. The homomorphism
\[
\Phi_h:\operatorname{Pic}^0(X)\longrightarrow\operatorname{Alb}(X)
\ \ \textnormal{given by } \ 
L\longmapsto a_X\bigl(h^{n-1}c_1(L)\bigr),
\]
is an isogeny. Consequently, if $L$ is non-torsion, then
\(
h^{n-1}c_1(L)\neq0
\)
in
\(
\operatorname{CH}^n(X)_{\mathbb{C}},
\)
and therefore $h\,c_1(L)\neq0$ in $\operatorname{CH}^2(X)_{\mathbb{C}}$.
\end{lemma}

\begin{proof}
Choose a smooth complete-intersection curve
$j\colon C\hookrightarrow X$ of class $h^{n-1}$, and put
$r=j^*\colon \operatorname{Pic}^0(X)\to\operatorname{Pic}^0(C)$.
Weak Lefschetz implies that $r$ has finite kernel.
By the projection formula and Picard--Albanese duality,
$\Phi_h$ identifies, up to sign, with
$r^\vee\circ\xi\circ r$, where $\xi$ is the canonical
principal polarization of $\operatorname{Pic}^0(C)$.
This is the homomorphism associated with
$r^*\mathcal O(\Theta)$, which is ample because $r$ is
finite onto its image. Thus $\Phi_h$ is an isogeny.

If $L$ is non-torsion, then $\Phi_h(L)$ is non-torsion,
so its image in
$\operatorname{Alb}(X)(\mathbb C)\otimes_{\mathbb Z}\mathbb C$
is nonzero. Applying the complexified Albanese homomorphism
therefore gives $h^{n-1}c_1(L)\neq0$.
The second assertion follows from
\(
h^{n-1}c_1(L)=h^{n-2}\bigl(h\,c_1(L)\bigr).
\)
\end{proof}

\begin{corollary}
\label{cor:irregularity-closed-rational}
The following conditions are equivalent:
\begin{enumerate}
\item $q(X)=h^1(X,\mathcal{O}_X)=0$;

\item $\Delta(\mathcal{F}_\eta)=0$ for every nonzero
closed rational one-form $\eta$ on $X$.

\end{enumerate}
\end{corollary}

\begin{proof}
If $q(X)=0$, then $\operatorname{Pic}^0(X)=0$.
Theorem~\ref{thm:closed-rational-discrepancy} therefore
implies $\Delta(\mathcal{F}_\eta)=0$ for every nonzero
closed rational one-form $\eta$.
Thus $(1)$ implies $(2)$.

Suppose that $(2)$ holds.
If $q(X)>0$, the positive-dimensional complex abelian
variety $\operatorname{Pic}^0(X)$ contains a non-torsion
point $L$.
Choose a very ample class $h$.
Lemma~\ref{lem:ample-albanese-isogeny} gives
\(
0\neq h\,c_1(L)
\in
\operatorname{Pic}^0(X)_{\mathbb{C}}
\cdot\operatorname{CH}^1(X)_{\mathbb{C}}.
\)
This contradicts Corollary~\ref{cor:discrepancy-span}.
Hence $q(X)=0$.
\end{proof}

\section{Final remarks}

Let $X$ be a smooth connected projective variety over $\mathbb C$.
For a codimension-one foliation $\F$ on $X$, the Baum--Bott theorem
implies that $\Delta(\F)$ is homologically trivial. When $\F$ is
defined by a closed rational one-form, Theorem~\ref{thm:closed-rational-discrepancy} gives the stronger
conclusion
\[
   \Delta(\F)\in
   \operatorname{Pic}^0(X)_{\mathbb C}\cdot
   \mathrm{CH}^1(X)_{\mathbb C}
   \subseteq \mathrm{CH}^2(X)_{\mathrm{alg},\mathbb C}.
\]
Example~\ref{ex:elliptic} shows that the discrepancy can be nonzero modulo
rational equivalence even for a logarithmic one-form on a surface.
This leads us to ask whether algebraic triviality persists for
arbitrary codimension-one foliations.

\begin{conjecture}\label{conj:algebraic-discrepancy}
For every codimension-one foliation $\F$ on $X$, the discrepancy
$\Delta(\F)$ is algebraically equivalent to zero.
\end{conjecture}

On a surface the conjecture follows directly from the Baum--Bott
theorem, since every zero-cycle of degree zero on a smooth connected
projective surface is algebraically equivalent to zero. Thus the
first case in which the conjecture imposes a further condition is
that of foliations on threefolds.

Even for regular foliations, Conjecture~\ref{conj:algebraic-discrepancy}
would imply the rank-two case of the weak Bloch conjecture for flat
vector bundles; see \cite[Conjecture~1.2]{Langer2024}. Indeed, let $E$
be a rank-two vector bundle with a flat algebraic connection on a
smooth connected projective variety $B$. Write
$\pi\colon Y=\mathbb P_B(E)\to B$ for the projective bundle of
one-dimensional quotients and put $\xi=c_1(\mathcal O_Y(1))$.
The connection induces a regular horizontal foliation $\mathcal H$
of codimension one on $Y$; see \cite[\S\,1.2]{LorayPereiraTouzet2016}. Its normal bundle is $T_{Y/B}$, and its
Baum--Bott cycle is zero because $\mathcal H$ is regular. Hence the
relative Euler sequence and the projective bundle relation give
\[
   \Delta(\mathcal H)
   =c_1(T_{Y/B})^2
   =\bigl(2\xi-\pi^*c_1(E)\bigr)^2
   =\pi^*\bigl(c_1(E)^2-4c_2(E)\bigr).
\]
If Conjecture~\ref{conj:algebraic-discrepancy} holds, intersecting
with $\xi$ and applying $\pi_*$ shows that $c_1(E)^2-4c_2(E)$ is
algebraically trivial with complex coefficients, since
$\pi_*(\xi\cdot\pi^*\gamma)=\gamma$. The same conclusion holds with
rational coefficients by extension of scalars. Moreover, $\det E$
is flat, so its first Chern class is torsion in the N\'eron--Severi
group of $B$. It follows that
\[
   c_2(E)\in\mathrm{CH}^2(B)_{\mathrm{alg},\mathbb Q},
\]
as predicted by the weak Bloch conjecture.

Our second question concerns the residue on a connected component
of the singular locus. Let $\F$ have codimension $q$, let $\varphi$
be a homogeneous characteristic polynomial of degree
$q<d\leq n=\dim X$, and put $m=n-d$.

\begin{conjecture}\label{conj:rational-residue}
Let $S_\alpha$ be a connected component of $\operatorname{Sing}(\F)$.
If $\operatorname{Res}_{\varphi}(\F,S_\alpha)$ belongs to
$H_{2m}(S_\alpha,\mathbb Q)$, then it is algebraic.
\end{conjecture}

In~\cite[Theorem~D]{CorreaSeadeEnvelopes2026}, Corr\^ea and Seade
verify the rationality hypothesis in further cases under assumptions
on the transverse type of the foliation along a smooth compact
singular component.
The Nash comparison gives a reason to expect algebraicity.
For a polynomial with rational coefficients, Corollary~\ref{cor:bb-nash-algebraicity} expresses
the difference between the Baum--Bott residue and the pushforward
of the Nash residue as the class of a rational algebraic cycle.
Thus, in the remaining cases, it would suffice to establish
algebraicity for the corresponding Nash residue. The question is
particularly relevant in higher codimension, where the pure powers
of $c_1$ no longer exhaust the monomials outside $I_q$.


\end{document}